\documentclass[leqno]{amsart}

\usepackage{amssymb}
\usepackage{amsmath}
\usepackage{amsthm}
\usepackage[latin1]{inputenc}
\usepackage{graphicx}
\usepackage{color}

\newtheorem{thm}{Theorem}[section]
\newtheorem{lem}[thm]{Lemma}

\newtheorem{prop}[thm]{Proposition}

\numberwithin{equation}{section}

\def\R{{\mathbb{R}}}

\newcommand{\eps}{\varepsilon}

\newcommand{\pd}[2]{\frac{\partial #1}{\partial #2}}

\def\p#1,#2,#3 {\rlap{\kern#1pt\raise#2pt\hbox{#3}}}

\title[Existence of flame balls]{Existence of radial stationary solutions for a system in combustion theory}

\author{J\'er\^ome Coville} 

\address{ INRA, Equipe BIOSP\\
Centre de Recherche d'Avignon\\ 
Domaine Saint Paul, Site Agroparc\\
84914 Avignon cedex 9, France
}

\email{jerome.coville@avignon.inra.fr}

\author{Juan D\'avila}
\address{J.~D\'avila: Departamento de Ingenier\'ia Matem\'atica and Centro de Modelamiento Matem\'atico\\
UMI 2807 CNRS\\
Universidad de Chile\\
Blanco Encalada 2120 - 5 Piso\\
Santiago - Chile}
\email{jdavila@dim.uchile.cl}

\date{January 26, 2011}

\begin{document}

\begin{abstract}
In this paper, we construct radially symmetric solutions of a nonlinear noncooperative elliptic system derived from a model for flame balls with radiation losses. This model is based on a one step kinetic reaction and our system is obtained by approximating the standard Arrehnius law  by an ignition nonlinearity, and by simplifying the term that models radiation. We prove the existence of 2 solutions  using degree theory. 

\end{abstract}

\maketitle

\section{Introduction}
This paper deals with radial solutions of the system of equations
\begin{align}
\label{system00}
\left\{
\begin{aligned}
\Delta u - \eps  u &= - v f(u) \quad \hbox{in $\R^3$}
\\
\Delta v &=  v f(u)  \quad \hbox{in $\R^3$}
\\
u & \ge 0 , \quad v \ge 0 \quad  \hbox{in $\R^3$}
\\
\lim_{|x|\to\infty} u(x) & = 0  , \quad 
\lim_{|x|\to\infty} v(x) = 1  ,
\end{aligned}
\right.
\end{align}
where $\eps \ge 0$ and $f$ is an ignition type nonlinearity, that is,  there exists  $0<\theta<1$ such that
\begin{align}
\label{h1}
f(t) = 0 \quad \forall t \le \theta, \quad  f(t) > 0 \quad \forall t > \theta.
\end{align}

System \eqref{system00} arises as a model problem for some reaction diffusion systems in combustion theory. 
We will describe the connection with these models at the end of the section.

We call $(u,v) = (0,1)$ the trivial solution of \eqref{system00}. 
Note that for any solution  $(u,v)$ of \eqref{system00} one has
$$
0 \le u \le 1 \quad \hbox{and} \quad 0 \le v \le 1,
$$
as can be deduced from the maximum principle.

We are mainly interested in existence and multiplicity of nontrivial radial solutions of \eqref{system00}. 
The first observation in this direction is that for $\eps>0$ large \eqref{system00} has no nontrivial solutions.
Indeed, by the maximum principle $u + v \le 1$ and then
$$
-\Delta u \le (1-u) f(u) - \eps u \quad \hbox{in $\R^3$.}
$$
If $\eps>0$ is large then $ (1-u) f(u) - \eps u \le 0 $ for all $u \in [0,1]$ and we conclude that $u \equiv 0$. Better estimates for
the quantity
$$
\eps^* = \sup \{ \eps>0 : \hbox{\eqref{system00} has a nontrivial radial solution} \} 
$$
are given in Section~\ref{estimate on eps*}.



\bigskip

Suppose that $u,v$ is a nontrivial radial solution.  
Then for some $r>0$ we must have $u(r) >\theta$. Otherwise $f(u) \equiv 0$ and then $\Delta v=0$, and 
by the Liouville theorem $v\equiv 1$. Since $0\le u+v\le 1$ this would imply  $u\equiv 0$ and then $u,v$ is trivial.
So for a nontrivial solution $u,v$, since $\lim_{r\to\infty} u(r) = 0$, there is $\beta_u>0$ such that 
$$
u(\beta_u) = \theta \quad\hbox{and} \quad u(r) <\theta \quad \forall r>\beta_u.
$$
The parameter $\beta_u$, which we will write simply as $\beta$, then serves to distinguish different solutions. 

\begin{thm} \label{thm existence1}
Assume $f : [0,1] \to [0,+\infty)$ is continuous, satisfies \eqref{h1} and for some $C>0$ 
\begin{align}
\label{h2}
 f(t) \le C f(u) \quad \forall  t \le u , \quad t, u \in [\theta,1] .
\end{align}
Then $\eps^*>0$ and there exist $0<\eps_1<\eps^*$ such that
for $0<\eps<\eps_1$ there are at least 2 solutions of \eqref{system00}. One of them has bounded $\beta$ as $\eps \to 0$ and the other has $\beta$ in the range  $\delta/\sqrt \eps \le \beta \le 1/(\theta\sqrt\eps)$ as $\eps\to 0$, where $\delta>0$ is fixed.
\end{thm}

A very natural and interesting question is the stability of the solutions constructed in Theorem~\ref{thm existence1}. Based on the works \cite{Buckmaster Joulin Ronney 1990, Buckmaster Joulin Ronney 1991}   we conjecture that the solution with bounded $\beta$ is unstable and that in part of branch of solutions with large $\beta$ the solution is stable, at least with respect to radial perturbations.

In some cases one may want to consider a discontinuous nonlinearity, such as the Heavisde function $f(u) = \chi_{[u>\theta]}$. 
With this example in mind we introduce the following hypothesis
\begin{align}
\label{h1b}
\hbox{$f$ is continuous in $(\theta,1]$ and $\lim_{u\to \theta+} f(u)$ exists} .
\end{align}


\begin{thm} 
\label{thm existence2}
Assume $f : [0,1] \to [0,+\infty)$ satisfies \eqref{h1}, \eqref{h2}, and \eqref{h1b}.
Then $\eps^*>0$ and there exist $0<\eps_1<\eps^*$ such that
for $0<\eps<\eps_1$ there are at least 2 solutions of \eqref{system00}. One of them has bounded $\beta$ as $\eps \to 0$ and the other has $\beta$ in the range  $\delta/\sqrt \eps \le \beta \le 1/(\theta\sqrt\eps)$ as $\eps\to 0$, where $\delta>0$ is fixed.
\end{thm}

The solution constructed in 
Theorem~\ref{thm existence2} is such that the set $\{r \in [0,\infty) : u(r) = \theta \}$ is finite, and hence the equation \eqref{system00} holds a.e.\@ in $\R^3$.

The motivation to consider a discontinuous nonlinearity is
only mathematical. However the example  $f(u)= \chi_{[u > \theta]}$ is interesting since it provides a situation where explicit calculations are possible. Theorem~\ref{thm existence2} shows that in part the conclusions obtained for  $f(u)= \chi_{[u > \theta]}$ remain valid for more general non-linearities.

\begin{figure}
\begin{center}
\label{fig1}
\hbox{
\includegraphics[width=10cm]{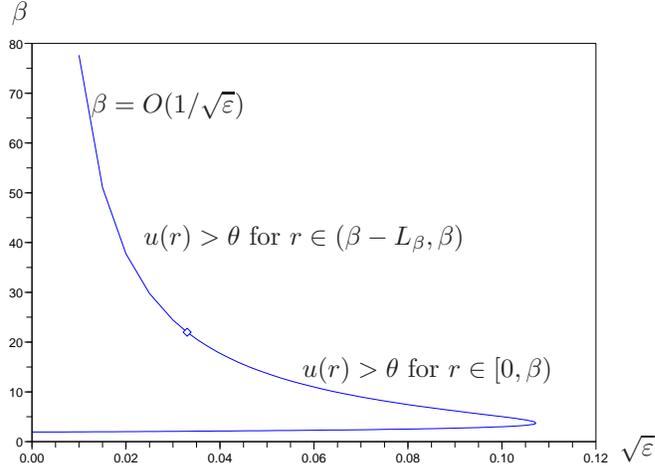}
\p -30,20,{$\sqrt\eps$}
\p -260,185,{$\beta$}
\p -150,50,{$u(r)>\theta$ for $ r \in [0,\beta)$}
\p -210,100,{$u(r)>\theta$ for $r \in (\beta-L_\beta,\beta)$}
\p -230,150,{$\beta =O(1/\sqrt\eps)$}
}
\caption{Bifurcation diagram for $f= \chi_{[u > \theta]}$ with $\theta=0.5$; $\beta$ is in the vertical axis and $\sqrt\eps$ in the horizontal axis}
\end{center}
\end{figure}

For $f(u)= \chi_{[u > \theta]}$ explicit calculations lead to an  equation for $\beta$ and $\eps$ in order for a radial solution to exist. In Figure~\ref{fig1} we show the numerical solution for this relation when $\theta=0.5$, with $\beta$ in the vertical axis and $\sqrt\eps$ in the horizontal axis.  It shows that for $0<\eps<\eps^*$ there are 2 solutions. Solutions in the lower branch satisfy $u>\theta$ in $[0,\beta)$ that is, the reaction takes place in the ball of radius $\beta$. The same happens for points in the upper branch which are to the right of the special point marked in the graph. To the left of that point the solution satisfies $u>\theta$ in an annulus of the form $ r\in (\beta - L_\beta ,\beta)$. Thanks to the explicit form of the relation bewteen $\beta$ and $\eps$ we can compute the asymptotic behavior of the curve as $\beta\to\infty$, and we find that
$$
\lim_{\beta \to \infty} \sqrt\eps \beta = a_0 >0  \quad \hbox{and} \quad \lim_{\beta\to\infty} L_\beta = L_0 >0,
$$
where $a_0,L_0$ is the unique solution of the system of equations
\begin{align*}
\theta \left( \frac{a}{\tanh(a)} -1 \right) + \theta( 1 + a ) &= 1
\\
1 - \frac{1}{L \tanh(L)} + \frac{1}{L \sinh(L)} & = \theta (1 + a) .
\end{align*}
%
%

Because of the information on the Heaviside nonlinearity one can conjecture  that for general $f$ there should be a similar relation for $\beta$ and $\eps$ as $\beta \to+\infty$.
We present in Section~\ref{Apriori estimates} 
nonexistence results for general ignition nonlinearities satisfying \eqref{h1} and \eqref{h2}, that capture this relation, and
roughly speaking say that no solution can exist if $\sqrt \eps \beta$ is either too large or too small, provided $\beta$ is taken large enough.
Using these nonexistence results and degree theory we can give a proof of Theorem~\ref{thm existence1}. This is done in Section~\ref{Existence of solutions}. In Section~\ref{sect exsit discont} we give the proof of Theorem~\ref{thm existence2} by approximating the discontinuous nonlinearity by continuous ones. 
Section~\ref{The Heaviside ignition} is devoted to the explicit computations for the Heaviside function.
Finally Section~\ref{estimate on eps*} contains a finer estimate of $\eps^*$.

As mentioned before, system \eqref{system00} arises in connection with some models in combustion theory, more precisely, in the flame ball problem for a weakly premixed gas sensitive to radiative heat losses. In such a mixture, it is known that apparently stationary spherical structures appear, which are called flame balls \cite{Lewis and Elbe, Ronney, Buckmaster Joulin Ronney 1990, Buckmaster Joulin Ronney 1991}.

In \cite{Lewis and Elbe, Ronney, Buckmaster Joulin Ronney 1990, Buckmaster Joulin Ronney 1991, Kagan Sivashinsky 1997, Kagan Minaev Sivashinsky 2004, Minaev Kagan Joulin Sivashinsky 2001} the following reaction diffusion system has been used to model a combustion process where flame balls can appear:
\begin{equation}
\left \{
\begin{aligned}
\rho C_p\frac{\partial T}{\partial t} &= \nabla \cdot (\lambda\nabla T) +\frac{QBY\rho}{m}e^{-\frac{E}{RT}} - q(T)\\
\rho\frac{\partial Y}{\partial t} &= \nabla \cdot (\mu\nabla Y) -BY\rho e^{-\frac{E}{RT}}
\\
Y&\to Y_\infty, \quad T\to T_\infty\quad \text{ as }\quad |x|\to \infty
\end{aligned} 
\right .\label{bjr}
\end{equation}
where $T$ is the temperature, $Y$ the reactant concentration, $Y_\infty>0$, $T_\infty >0$ are the reactant concentration and tempereature at infinity, and $C_p$, $R$, $Q$ and $m$ are respectively the specific heat capacity at constant pressure, the perfect gas constant, the chemical heat release and the molecular mass of the reactant. The term $q(T)$ represents radiative losses. The reaction is characterized by the one-step Arrhenius kinetics $\sim Be^{\frac{-E}{RT}}$ where $B$ is a constant. Furthermore, the hydrodynamics effects are neglected, i.e the density $\rho$, the thermal conductivity $\lambda$ and  the diffusion coefficient $\mu$ are constant. See also numerical simulations in \cite{Bockhorn Frochlich Schneider 1999, Kagan Sivashinsky 1997, Kagan Minaev Sivashinsky 2004}. 

After the seminal work \cite{zeldovich}, the traveling front problem for systems like \eqref{bjr} has been  investigated by several authors, for example, 
\cite{giovangigli,Glangetas-Roquejoffre,roques}.

In absence of radiation, i.e.\@ when $q=0$, there are many works dealing with  \eqref{bjr}, see for instance
\cite{berestycki-larrouturou,Berestycki-Nicolaenko-Scheurer,Lederman-Roquejoffre-Wolanski,marion,sagon} and the references therein.
Also we remark that the stationary version of \eqref{bjr} without radiation, leads to 
system \eqref{system00} with $\eps=0$ which reduces to a scalar equation, since $v=1-u$. There is a huge amount of literature concerning existence of radial ground states for semilinear equations, so we mention here only some classical references  \cite{berestycki-lions-1,kwong} on the problem in entire space.
When the problem is treated in a bounded domain see the book \cite{Aris} and \cite{Du-Lou,Kapila-Matkowsky-Vega} for multiplicity results in the case of Arrehnius non-linearity.
The paper \cite{Brindley-Jivraj-Merkin-Scott} contains interesting numerical computations of the bifurcation diagram in the case of the full coupled system in a interior of a sphere.

A common simplification of \eqref{bjr} under the assumption of large activation energy, that is $E>>1$, is to assume that the source term for the reaction is concentrated on a very thin layer, typically a sphere. This approach is taken for example in \cite{Buckmaster Joulin Ronney 1990, Buckmaster Joulin Ronney 1991,  Lee Buckmaster 1991, Minaev Kagan Joulin Sivashinsky 2001}
and leads to the free boundary problem
 \begin{equation}
  \left \{
   \begin{aligned}
    \frac{\partial T}{\partial t} & = \Delta T +Be^{-\frac{E}{2T^*}}\delta(r-R(t)) - q(T) \quad \text{in }\quad \R^3\\
    \frac{\partial Y}{\partial t} & = \frac{1}{Le}\Delta Y -Be^{-\frac{E}{2T^*}}\delta(r-R(t)) \quad  \text{in }\quad \R^3\\
      Y & \equiv 0\quad \text{ in }\quad B(0,R(t))\\
      q & \equiv 0\quad \text{ in }\quad \R^3\setminus B(0,R(t))\\
Y&\to Y_\infty, \quad T\to T_\infty\quad \text{ as }\quad |x|\to \infty
   \end{aligned} 
  \right .\label{bjr-fb}
\end{equation}
where $R(t)$ is the radius of  the front where the reaction takes place, $\delta$ is the Dirac measure and $T^*$ is the front temperature. In \cite{Buckmaster Joulin Ronney 1990, Buckmaster Joulin Ronney 1991,Buckmaster Lee 1991} the authors analyze the stability of stationary solutions of \eqref{bjr-fb}. In a similar framework, existence and stability of flame balls and travelling flame balls have also been studied in \cite{Minaev Kagan Joulin Sivashinsky 2001,Guyonne Lorenzi 2007,Guyonne Noble 2007,Shah:2000qf,Berg:2006qy}.


We arrive at \eqref{system00}
by introducing the following simplifications:
\begin{itemize}
\item[a)]
assume the radiative loss to be linear, i.e.,  $q(T)= a (T-T_\infty)$ where $a>0$,
\item[b)]
approximate $ e^{-\frac{E}{RT}}$ by  $ \eta(T-T_0)BYe^{-\frac{E}{RT}}$ where  $T_0$ is an activation temperature and $\eta$ is a cut-off function satisfying $\eta>0$ in $\R^+$ and $\eta\equiv 0$ in $\R^-$.
\end{itemize}
As in \cite{Kagan Minaev Sivashinsky 2004,Shah:2000qf} one can model radiative heat losses using Stefan's law $q(T) = \eps (T^4 - T_\infty^4)$ for some constant $\eps>0$. When $T$ is close to $T_\infty$ we can write $q(T) \approx 4 \eps T_\infty^3 (T-T_\infty)$.
As a step towars understanding more general situations, we assume that this linear relation holds for all $T$, that is, we assume a). 
Other linear or piecewise linear approximations have been used before, for instance in \cite{Brailovski-Sivashinsky,Minaev Kagan Joulin Sivashinsky 2001,Shah:2000qf}.
Assumption b) corresponds the a standard approximation in combustion theory to avoid the {\em cold boundary difficulty}, see \cite{Berestycki-Larrouturou-Roquejoffre}.

After introducing dimensionless variables $u$, $v$, corresponding to temperature and reactant concentration, the stationary version of \eqref{bjr} becomes
\begin{equation}
\label{system02}
\left \{ 
\begin{aligned}
& \Delta u + v g(u)- c u  = 0 \quad \hbox{in }  \R^{3}\\
& \frac{1}{Le}\Delta v - v g(u) = 0  \quad \hbox{in  } \R^{3}\\
& u \to 0 , \quad v \to v_\infty  \quad \hbox{as  }|x|\to +\infty
\end{aligned}
\right.
\end{equation}
where $Le>0$ is the Lewis number, $c>0$, $v_\infty>0$ and $g$ is an ignition type function, that is, there is $\theta>0$ such that   $g(u)=0$  if $u \le \theta$.

We stress that our results are valid for any value of $Le>0$. Indeed, since we are considering stationary solutions, the following change of variables will allow us to assume that $Le=1$. Letting $v(x) = v_\infty \, \tilde v(\sqrt{Le} \,  x) $, $u(x) =  \frac{v_\infty}{Le} \tilde u(\sqrt{Le} \,  x) $ transforms system
\eqref{system02} into \eqref{system00}, where $\eps = \frac{c}{Le}>0$ and $f(u) = g( \frac{v_\infty}{Le}u)$. We observe that $f$ is still is an ignition type nonlinearity.

\section{Apriori estimates}
\label{Apriori estimates}
The purpose in this section is to establish nonexistence results in some ranges of the parameters. 

Given a nontrivial solution $u,v$ of \eqref{system00} let $\beta>0$ be such that
$$
u(\beta) = \theta \quad\hbox{and} \quad u(r) <\theta \quad \forall r>\beta.
$$
Setting $\tilde u(r) = u(\beta r)$, $\tilde v(r) = v(\beta r)$ these new functions satisfy

\begin{align}{
\label{system01}\left\{
\begin{aligned}
&\frac{\partial}{\partial r}\left(r^2\frac{\partial u}{\partial r} \right) - \eps r^2 \beta^2 u = - r^2\beta^2 v f(u) &\quad \hbox{in $\R^+$}
\\
&\frac{\partial}{\partial r}\left(r^2\frac{\partial v}{\partial r} \right) =  r^2\beta^2 v f(u) &\quad \hbox{in $\R^+$}
\\
&\lim_{r\to\infty} u(r)  = 0  , \qquad 
\lim_{r\to\infty} v(r) = 1 
\end{aligned}
\right.}
\end{align}
In the sequel we will study \eqref{system01} in the following set of functions
$$
\mathcal S = \{ ( u,v) : \hbox{$u,v$ are $C^1([0,\infty))$ and $u(1)  = \theta$ and $u(r)<\theta$ for $r>1$}\}.
$$

Let $h_0:\R\to\R$ be such that
\begin{align}
\label{f0 1}
& h_0(t) = 0 \quad \forall t \le \theta \quad \hbox{and} \quad  h_0(t) > 0 \quad \forall t > \theta
\\
\label{f0 2}
& \hbox{$h_0$ is nondecreasing.}
\end{align}
We consider now functions $f:\R\to\R$ such that
\begin{align}
\label{f1}
& f(t) = 0 \quad \forall t \le \theta, 
\\
\label{f2}
& f \ge h_0 \quad \hbox{in $\R$}
\\
\label{f3}
& f(t) \le C_0 f(u) \quad \forall \theta \le t \le u \le 1.
\end{align}

\begin{lem}
\label{nonexistence large epsilon beta2}
There is no solution in $\mathcal S$ to \eqref{system01} if 
$$
\beta \sqrt \eps > \frac{1}{\theta}.
$$
\end{lem}
\proof[\bf Proof]
Since
$$
\Delta u - \eps \beta^2 u = 0 \quad \hbox{for $r>1$}, \quad u(1) = \theta,
\quad
\lim_{r\to + \infty} u(r) = 0,
$$
we have an explicit formula
$$
u(r) = \theta \frac{e^{-(r-1)\beta \sqrt \eps}}{r} \quad \hbox{for all $r\ge 1$.}
$$
>From this we find
\begin{align}
\label{u prime 1}
u'(1) = - \theta ( 1 + \beta \sqrt \eps).
\end{align}
Similarly, since $\Delta v = 0 $ for $r\ge 1$ and $\lim_{r\to\infty} v(r) =1 $ we have
$$
v(r) = 1 - {\frac{\gamma}{r}} \quad \forall r\ge1
$$
where $0 < \gamma < 1$. This yields
$$
v'(1) = \gamma \in (0,1).
$$
Integrating the equation for $v$ in (0,1) implies
\begin{align}
\label{v prime 1}
\beta^2 \int_0^1 s^2 v f(u) \, d s = v'(1) = \gamma.
\end{align}
Let $r_0 \in [0,1)$ be such that $u'(r_0)=0$. Then integrating the equation for $u$ in $(r_0,1)$ yields
$$
u'(1) = \eps \beta^2 \int_{r_0}^1 s^2 u \, d s - \beta^2 \int_{r_0}^1 s^2 v f(u) \, d s.
$$
This formula together with \eqref{u prime 1} and  \eqref{v prime 1} gives
$$ 
v'(1) \ge  \beta^2 \int_{r_0}^1 s^2 v f(u) \, d s =  \eps \beta^2 \int_{r_0}^1 s^2 u \, d s  +  \theta ( 1 + \beta \sqrt \eps) 
$$
and it follows that
$$
\theta \beta \sqrt \eps \le 1.
$$
\qed

\begin{lem}
\label{nonexistence small epsilon beta2}
There is $\beta_0>0$, $\delta>0$ depending only on $h_0$, $\theta$, $C_0$ such that for all $\beta \ge \beta_0$, $\beta \sqrt \eps \le \delta$ and all $f$ satisfying \eqref{f1}, \eqref{f2}, \eqref{f3} there is no solution in $\mathcal S$ to the system \eqref{system01}.
\end{lem}
\proof[\bf Proof]
We treat the case $\eps>0$ since the situation $\eps=0$ is similar. 

As before
$$
u(r) = \theta \frac{e^{-(r-1)\beta \sqrt \eps}}{r} \quad \hbox{for all $r\ge 1$.}
$$
>From this we find
\begin{align}
\label{u prime 1 2}
u'(1) = - \theta ( 1 + \beta \sqrt \eps).
\end{align}
Similarly
$$
v(r) = 1 - \frac{\gamma}{r} \quad \forall r\ge1
$$
where $0 < \gamma < 1$. This yields
$$
v'(1) = \gamma \in (0,1).
$$
Integrating the equation for $u$ we see that
\begin{align}
\label{001b}
r^{2} u'(r) = \eps \beta^2 \int_0^r s^2u \, d s  - \beta^2 \int_0^r s^{2} v f(u) \, d s .
\end{align}
>From this and \eqref{u prime 1 2} it follows that
\begin{align}
\label{004b}
\beta^2  \int_0^1 s^{2} v f(u) \, d s = \theta ( 1+ \beta \sqrt\eps)  + \eps \beta^2 \int_0^1 s^2 u \, d s \le \theta + \delta(\theta + \delta).
\end{align}
Integrating the equation for $v$ in (0,1) implies
\begin{align}
\label{v prime 1b}
\beta^2 \int_0^1 s^2 v f(u) \, d s = v'(1) = \gamma
\end{align}
and combined  with \eqref{004b} yields 
\begin{align}
\label{gamma}
\gamma \le \theta + \delta(\theta + \delta).
\end{align}

Integrating the equation for $v$ in $(0,r)$ we obtain
\begin{align}
\label{v prime r}
v'(r) = \beta^2 \int_0^r (s/r)^2  v(s) f(u(s)) \, d s.
\end{align}
Integrating this on $(0,R)$ with $0\le R \le 1$ yields
\begin{align*}
v(R) - v(0) & = \beta^2 \int_0^R  \int_0^r (s/r)^2  v(s) f(u(s)) \, d s \, d r 
\\
& =  \beta^2 \int_0^R  \int_s^R r^{-2} \, d r s^2  v(s) f(u(s)) \, d s 
\\
& =  \beta^2 \int_0^R  (s^{-1} - R^{-1}) s^2  v(s) f(u(s)) \, d s .
\end{align*}
In particular, with $R=1$
$$
\beta^2 \int_0^1  (s - s^2)  v(s) f(u(s)) \, d s = v(1) - v(0) \le 1
$$
This and \eqref{v prime 1b} give
$$
\beta^2 \int_0^1  s   v(s) f(u(s)) \, d s \le 2.
$$
Now going back to \eqref{v prime r} we obtain
\begin{align}
\label{derivative of of}
v'(r) & =  \beta^2 \int_0^r (s/r)^2  v(s) f(u(s)) \, d s \le   \beta^2 \int_0^r (s/r) v(s) f(u(s)) \, d s \le  \frac{2}{r } \quad 0< r \le 1.
\end{align}
Now consider the function $z = u+ v$ which satisfies $\Delta z = \eps \beta^2 u $ in $\R^3$. Integrating the relation $(s^2 z')' = \eps \beta^2 s^2 u$ in $(0,r)$ 
$$
z'(r) = \eps \beta^2 \int_0^r (s/r)^2 u(s) \,  d s 
$$
so that 
\begin{align}
\label{bound z prime}
0 \le z'(r) \le \delta^2 r \quad r\ge 0.
\end{align}
It follows that
\begin{align}
\label{derivative u}
-\frac{2}{r} \le u'(r) \le \delta^2 r \quad 0<r\le 1.
\end{align}
Observe that integrating \eqref{bound z prime} on $(r,1)$ we find
$$
0 \le z(1) - z(r) \le \delta^2
$$
and since $z(1) = \theta + 1 - \gamma$
$$
1 + \theta- \gamma - \delta^2 \le z(r) \le   1 + \theta - \gamma \quad 0 \le r \le 1.
$$
In particular, using \eqref{gamma},
\begin{align}
\label{v from below}
v = z- u \ge 1 + \theta  - \gamma - \delta^2 - u \ge 1 - \delta ( \theta + 2 \delta) - u \quad 0 \le r \le 1.
\end{align}

\medskip \noindent {\bf Step 1.} For any $r_0>0$ there exists $\beta_1(r_0)$ depending on $r_0$, $h_0$ and $\theta$ only such that 
\begin{align}
\label{006b}
\max_{[r_0,1]} u < \frac{1+\theta}{2} \quad \forall \beta \ge \beta_1 .
\end{align}
To prove this, suppose that $\max_{[r_0,1]} u \ge \frac{1+\theta}{2} $ and let $r \in [r_0,1]$ be such that $u(r) =   \frac{1+\theta}{2} $. Let $M = 2/r_0$ so that from \eqref{derivative u}
$$
|u'(s)| \le M \quad r_0 \le s \le 1.
$$
Then 
\begin{align}
\label{u bound lipschitz}
\frac{1+\theta}{2} - \frac{1-\theta}{4} \le  u(s) \le \frac{1+\theta}{2} + \frac{1-\theta}{4} \quad \forall  s \in [r_0,1] \quad \hbox{such that} \quad |s-r| \le \frac{1-\theta}{4 M}.
\end{align}
Using \eqref{v prime 1b} and since $f \ge h_0$ 
\begin{align*}
1 & \ge \beta^2 \int_{r_0}^1  s^2 v(s) f(u(s)) \, d s \ge  \beta^2 \int_{r_0}^1  s^2 v(s) h_0(u(s)) \, d s .
\end{align*}
If  $r\le \frac{1+r_0}{2}$ let $I = [ r - \frac{1-\theta}{4M}, r]$ and if  $r\ge \frac{1+r_0}{2}$ let $I = [ r, r +   \frac{1-\theta}{4M}]$. Using \eqref{v from below} and \eqref{u bound lipschitz}:
\begin{align*}
1 &  \ge 
\beta^2 \int_I  s^2 v(s) f(u(s)) \, d s 
\\
& \ge \beta^2 r_0^2 \frac{1-\theta}{4 M} h_0\left( \frac{1+\theta}{2} - \frac{1-\theta}{4}  \right)  \left( 1 - \delta(\theta+ 2\delta) -  \frac{3+\theta}{4}\right) .
\end{align*}
We take $\delta>0$ small only depending on $\theta$ such that 
\begin{align}
\label{choice delta}
1 - \delta( \theta+ 2\delta) -  \frac{3+\theta}{4} \ge \frac{1-\theta}{8} .
\end{align}
Then
$$
1  \ge \beta^2 r_0^3 \frac{( 1-\theta)^2}{64} h_0\left( \frac{1+\theta}{2} - \frac{1-\theta}{4}  \right) 
$$
and the claim follows.

\bigskip \noindent {\bf Step 2.} 
For any $r_0>0$ there exists $C>0$ such that
$$
(\max_{[r_0,1]} u -\theta) h_0\left(\frac{\max_{[r_0,1]} u-\theta}{2} + \theta \right)  \le \frac{C}{\beta^2}.
$$
The constant $C$ depends only on $r_0$ and $\theta$. The conclusion from this is that
\begin{align}
\label{u-theta to zero}
\limsup_{\beta\to\infty} \max_{[r_0,1]} ( u - \theta) \le 0,
\end{align}
and this is uniform with respect to $u$ and $f$.

\medskip
As before, setting  $M = 2/r_0$ we have  $ |u'(s)| \le M $ for $ r_0 \le s \le 1$ by \eqref{derivative u}. Let $r_m \in [r_0,1]$ be a point such that $u(r_m) =  \max_{[r_0,1]} u$ and let  us write  $u_m = u(r_m) = \max_{[r_0,1]} u$. Then
$$
u(r) \ge \frac{u_m-\theta}{2} + \theta \quad \forall r \in [r_0,1] \quad \hbox{such that} \quad  |r-r_m| \le \frac{u_m - \theta}{2M}.
$$
Since $f\ge h_0$, using \eqref{v prime 1b}, \eqref{v from below}, \eqref{006b}, \eqref{choice delta} we obtain
\begin{align*}
1 & \ge \beta^2  \int_0^1 s^{2} v(s) f(u(s)) \, d s \ge  \beta^2  \int_0^1 s^{2}  v(s) h_0(u(s)) \, d s
\\
& \ge \beta^2   r_0^2  \frac{u_m - \theta}{2M} h_0\left(\frac{u_m-\theta}{2} + \theta \right)   \left( 1 - \delta(\theta+ 2 \delta) -  \frac{1+\theta}{2}\right)
\\
& \ge \beta^2   r_0^2  \frac{u_m - \theta}{2M} h_0\left(\frac{u_m-\theta}{2} + \theta \right)  \frac{1-\theta}{8}.
\end{align*}
Thus there exists $C$ depending only on $r_0$, $\theta$ such that for $\beta \ge \beta_1$
$$
(u_m-\theta) h_0\left(\frac{u_m-\theta}{2} + \theta \right)  \le \frac{C}{\beta^2}.
$$
This proves the claim.

\bigskip \noindent {\bf Step 3.} For $0\le a \le 1$
\begin{align}
\label{integral gradient}
\int_a^1 r^2 (u')^2 \, d r \le 2 a + \delta^2 +  \max_{[a,1]} (u - \theta).
\end{align}
Indeed multiplying the equation 
\begin{align}
\label{eq u by r2}
 - (r^2 u')' = -  \eps \beta^2 r^2 u + \beta^2 r^2 v f(u)
\end{align}
by $(u-\theta)$ and integrating over $(a,1)$ we get
$$
\int_a^1 r^2 (u')^2 \, d r =  r^2 u' (u-\theta)\Big|_a^1
- \eps \beta^2 \int_a^1 r^2 u (u-\theta) \, d r + \beta^2 \int_a^1 r^2 v f(u) (u-\theta) \, d r.
$$
By \eqref{derivative u} 
$$
r^2 u' (u-\theta)\Big|_a^1 = - a^2 u'(a)( u(a) - \theta) \le 2 a.
$$
Since $u \le 1$ we have
$$
\eps \beta^2  \int_a^1 r^2 u (u-\theta) \, d r  \le \delta^2
$$
Finally using \eqref{v prime 1b}
$$
\beta^2 \int_a^1 r^2 v f(u) (u-\theta) \, d r \le \max_{[a,1]} (u - \theta) \beta^2 \int_a^1 r^2 v f(u)  \, d r \le \max_{[a,1]} (u - \theta) .
$$

\bigskip \noindent {\bf Step 4.} 
We finish the proof of the Lemma using a modification of Pohozaev's identity.  Let $0<a<1$ be a fixed number.
Multiplying \eqref{eq u by r2} by $(r-a) u'$ and integrating over $(a,1)$ 
\begin{align}
\label{pohozaev}
- \int_a^1 (r^2 u')' (r-a) u' \, d r = \eps \beta^2 \int_a^1 r^2 (r-a) u u'\, d r  +  \beta^2 \int_a^1 r^2 (r-a) v f(u) u' \, d r.
\end{align}
A computation shows that the left hand side is given by
\begin{align}
\nonumber
\int_a^1 (r^2 u')' (r-a) u' \, d r & = \frac{1}{2} r^2 (r-a) (u')^2 \Big|_a^1 + \frac{1}{2} \int_a^1 r^2 (u')^2 \, d r - a \int_a^1 r (u')^2 \, d r
\\
\label{pohozaev left}
&=
 \frac{1}{2} (1-a)  u'(1)^2 + \frac{1}{2} \int_a^1 r^2 (u')^2 \, d r - a \int_a^1 r (u')^2 \, d r.
\end{align}
Thanks to \eqref{integral gradient} we find
$$
a \int_a^1 r (u')^2 \, d r \le  \int_a^1 r^2 (u')^2 \, d r \le 
2 a + \delta^2 +  \max_{[a,1]} (u - \theta).
$$
Since $u'(1) = - \theta ( 1+ \sqrt \eps \beta) $ we obtain from \eqref{pohozaev left}, \eqref{integral gradient} and the previous estimate
\begin{align}
\label{b0020}
\left| \int_a^1 (r^2 u')' (r-a) u' \, d r - \frac{1}{2}(1-a) \theta^2 ( 1+ \sqrt \eps \beta)^2 \right| \le 4 a + 2 \delta^2 +  2 \max_{[a,1]} (u - \theta).
\end{align}
Let us compute
$$
 \eps \beta^2 \int_a^1 r^2 (r-a) u u'\, d r = \frac{1}{2} \eps \beta^2 r^2 ( r-a) u^2\Big|_a^1 - \frac{1}{2} \eps \beta^2  \int_a^1 ( 3 r^2 - 2 a r) u^2 \, d r.
$$
Recalling that $\sqrt\eps \beta \le \delta$ this shows that
\begin{align}
\label{b001}
\left|  \eps \beta^2 \int_a^1 r^2 (r-a) u u'\, d r \right| \le 3 \delta^2.
\end{align}
Define
$$
F(u) = \int_\theta^u f(t) \, d t.
$$
Then
\begin{align*}
\beta^2 \int_a^1 r^2 (r -a)  v f(u) u' \, d r 
& = \beta^2 r^2 (r-a) v F(u) \Big|_a^1 - \beta^2 \int_a^1 (3r^2 - 2 a r) v  F(u) \, d r 
\\
& \qquad - \beta^2 \int_a^1 (r^3 - a r^2) v'  F(u) \, d r
\\
& =  \beta^2 \int_a^1 (3r^2 - 2 a r) v  F(u) \, d r  - \beta^2 \int_a^1 (r^3 - a r^2) v'  F(u) \, d r . 
\end{align*}
These terms can be estimated using \eqref{f3}:
$$
F(u) \le C_0 f(u) ( u- \theta) \quad \forall \theta \le u \le 1
$$
and hence by \eqref{v prime 1b}
$$
 \beta^2 \int_a^1 r^2  v  F(u) \, d r \le C_0 \int_a^1 r^2 v f(u) (u-\theta) \, d r
 \le C_0 \max_{[a,1]} (u-\theta) 
$$
Similarly
$$
a \beta^2 \int_a^1 r  v  F(u) \, d r \le  \beta^2 \int_a^1 r^2  v  F(u) \, d r 
 \le C_0 \max_{[a,1]} (u-\theta) .
$$
To estimate the remaining terms we use \eqref{derivative of of}
$$
\beta^2 \int_a^1 r^3  v'  F(u) \, d r \le 2 \beta^2 \int_a^1 r^2   F(u) \, d r \le
\frac{16}{1-\theta} \int_a^1 r^2  v F(u) \, d r
$$
where the last inequality is consequence of \eqref{v from below}, \eqref{006b} and \eqref{choice delta} provided $\beta \ge \beta_1(a)$ ($a \in (0,1)$ is fixed).
Hence 
$$
\beta^2 \int_a^1 r^3  v'  F(u) \, d r  \le   \frac{16 C_0}{1-\theta}  \max_{[a,1]} (u-\theta),
$$
and this implies also
$$
 \beta^2 \int_a^1 a r^2 v'  F(u) \, d r \le \beta^2 \int_a^1 r^3  v'  F(u) \, d r  \le   \frac{16 C_0}{1-\theta}  \max_{[a,1]} (u-\theta).
$$
Therefore
\begin{align}
\label{b002}
\left| \beta^2 \int_a^1 r^2 (r -a)  v f(u) u' \, d r \right| \le K_1   \max_{[a,1]} (u-\theta) 
\end{align}
where $K_1 = 5C_0 + \frac{32 C_0}{1-\theta}$.

>From \eqref{pohozaev}, \eqref{b0020}, \eqref{b001},  \eqref{b002} we arrive at
\begin{align}
\label{contra}
\frac{1}{2}(1-a) \theta^2 ( 1+ \sqrt\eps\beta)^2 \le 4a + 5 \delta^2 + (K_1+2) \max_{[a,1]} (u-\theta) 
\quad \forall \beta \ge \beta_1(a) .
\end{align}
We fix $\delta>0$ even smaller if necessary so that
$$
5 \delta^2 \le \frac{1}{3}\frac{\theta^2}{4} .
$$
Then fix $0<a<\frac{1}{2}$ such that
$$
4 a \le \frac{1}{3}\frac{\theta^2}{4}.
$$
Then for \eqref{contra} yields
$$
\frac{1}{3}\frac{\theta^2}{4} \le  (K_1+2) \max_{[a,1]} (u-\theta) 
\quad \forall \beta \ge \beta_1(a) 
$$
which is not possible for $\beta$ large enough by \eqref{u-theta to zero}.
\qed


\section{Existence of solutions}
\label{Existence of solutions}

In this section we prove Theorem~\ref{thm existence1},
and throughout it we assume that $f$ is continuous and satisfies \eqref{h1} and \eqref{h2}.

We work in the Banach space
$$
X = \{ (u,\beta) : u  \in C^{1,\alpha}([0,1]) , \;  u(1) = 0 , \; \beta\in \R \}
$$
endowed with its natural norm
$$
\| (u,\beta) \|_X = \|u\|_{C^{1,\alpha}([0,1])} + |\beta| ,
$$
where $0<\alpha<1$ is fixed.

Define for $t \in [0,1]$
$$
f_t(u) = t f(u) + (1-t) ( u-\theta)^+ \quad t \in [0,1],
$$
where $s^+ = \max(s,0)$.
Let us consider  \eqref{system01} with nonlinearity $f_t$ that is
\begin{align}
\label{system epsilon t}
\left\{
\begin{aligned}
\Delta u - \eps \beta^2 u &= - \beta^2 v f_t(u) \quad \hbox{in $\R^3$}
\\
\Delta v &=  \beta^2 v f_t(u) \quad \hbox{in $\R^3$}
\\
\lim_{|x|\to\infty} u(x) & = 0  , \qquad 
\lim_{|x|\to\infty} v(x) = 1 
\end{aligned}
\right.
\end{align}
To apply the non-existence results of the previous section we need to exhibit a function $h_0$ satisfying \eqref{f0 1} and \eqref{f0 2}.
For this purpose define
$$
g(t) = \inf \ \{ r\in [0,1]\ : \  f \ge t \hbox{ in } [r,1] \} \quad \hbox{for all } t \ge 0,
$$
and
$$
h(r) = \inf \ \{ t\ge 0 \ : \  g(t) \ge r \} \quad \hbox{for all } r \in [0,1].
$$
The following properties then follow from these definitions:
\renewcommand{\labelenumi}{\alph{enumi})}
\begin{enumerate}
\item
$g$ is strictly increasing, continuous from the left and satisfies $\lim_{t\to0^+} g(t) = \theta$.
\item
$h(r) = $ for $r\in [0,\theta]$, $h(r)>0$ for $r\in (\theta,1]$, $h$ is nondecreasing and continuous in $[0,1]$.
\item
$f\ge h$ in $[0,1]$.
\end{enumerate}
The function 
$$
h_0 = \min( h , (u-\theta)^+ )
$$
satisfies \eqref{f0 1} and \eqref{f0 2} and the nonlinearity $f_t$ satisfies \eqref{f1}, \eqref{f2}, \eqref{f3} for all $t\in [0,1]$.

Since any solution to \eqref{system epsilon t} is bounded above by 1 we may redefine $f_t(u)$ as a constant for $u \ge 1$. Thus we may assume that for some constant $M>0$ we have
$$
|  f_t(u) | \le M  \quad \forall u \in \R, \quad \forall t \in [0,1].
$$

We define a nonlinear map $T_{\eps,t}:X \to X$ as follows. Let $(\tilde u, \tilde \beta) \in X$. Then solve
\begin{align}
\label{T:sol v}
\Delta v = \tilde \beta^2 v f_t(\tilde u + \theta) \quad \hbox{in $B_1$}, \quad v(1) + v'(1) = 1 \quad \hbox{on $\partial B_1$}.
\end{align}
This problem has a unique solution which can be found for example by minimizing
$$
\frac{1}{2} \int_{B_1} \left( |\nabla v|^2 + \tilde \beta^2 f_t(\tilde u + \theta) v^2 \right)  + \frac{1}{2} \int_{\partial B_1} v^2 - \int_{\partial B_1} v 
$$
with $v \in H^1(B_1)$. Since $f_t(\tilde u + \theta) \in L^\infty(B_1)$ by standard elliptic regularity $v \in C^{1,\alpha}([0,1])$. Then find $\beta>0$ and $u$ such that 
\begin{align}
\left\{
\label{T:sol u}
\begin{aligned}
& - \Delta u + \eps \tilde \beta^2 u   =   \beta^2 v f_t(\tilde u + \theta) - \eps \tilde \beta^2 \theta + { \beta^2 - \tilde \beta^2 }  \quad \hbox{in $B_1$}
\\
& \quad u(1)=0 , \qquad 
u'(1) = - \theta ( 1 + \tilde \beta \sqrt \eps) 
\end{aligned}
\right.
\end{align}
admits a solution. This problem has a unique solution $(u,\beta)$ which furthermore is $C^{1,\alpha}([0,1])$ (this assertion is verified in the proof of Lemma~\ref{T continuous compact} below). We define
$$
T_{\eps,t}(\tilde u, \tilde \beta) = ( u , \beta)  .
$$
Observe that $(u,\beta)$ is a fixed point of $T_{\eps,t}$ if and only if $(u+\theta,v)$ is a solution of \eqref{system epsilon t}.

\begin{lem}
\label{T continuous compact}
Write $T(u,\beta,\eps,t) = T_{\eps,t}(u,\beta)$. Then
$$
T : X \times [0,\infty) \times [0,1] \to X
$$
is compact. Moreover
\begin{align}
\label{bound T}
\| T(u,\beta,\eps,t) \|_X \le C  ( 1 + \sinh(\sqrt \eps  \beta)  + \beta^2 ) \quad \forall (u,\beta) \in X, \forall \eps\ge 0, \forall t \in [0,1].
\end{align}
\end{lem}
\proof[\bf Proof.]
First we remark that the solution $v$ to \eqref{T:sol v} satisfies
$$
0 \le v \le 1 \quad \hbox{in $B_1$}
$$
because $0$ is a subsolution and $1$ is a supersolution. Given $a \ge 0$ let
$Z_a$ be the solution to
$$
-\Delta Z_a + a Z_a = 0 \quad \hbox{in $B_1$} , \quad Z_a(1) = 1.
$$ 
Note that $Z_a$ is explicit:
$$
Z_a(r) = \frac{\sinh( ar ) }{r \sinh(a)} \quad \forall 0< r \le 1
$$
and hence
\begin{align}
\label{bounds Za}
\frac{a}{\sinh(a)} \le Z_a \le 1.
\end{align}

Let $(u,\beta) = T(\tilde u,\tilde \beta, \eps,t)$. Then choosing $a =  \eps \tilde \beta^2$, multiplying \eqref{T:sol u} by $Z_a$ and integrating in $B_1$ we find
$$
-u'(1)  = \beta^2 \int_{B_1}( v f (\tilde u+ \theta) + 1) - \tilde \beta^2 \int_{B_1} ( \eps \theta + 1) Z_a .
$$
>From the boundary condition in \eqref{T:sol u} it follows that $\beta$ is explicitly given by
\begin{align}
\label{formula beta}
\beta^2 = \frac{\theta(1+\sqrt\eps\tilde\beta)  +  \tilde \beta^2 \int_{B_1} ( \eps \theta + 1) Z_a }{ \int_{B_1}( v f (\tilde u+ \theta) + 1) Z_a},
\end{align}
where  $a =  \eps \tilde \beta^2$ . By \eqref{bounds Za}
\begin{align}
\label{integral v f Z_a}
\int_{B_1}( v f (\tilde u+ \theta) + 1) Z_a \ge \frac{\sqrt \eps \tilde \beta}{\sinh(\sqrt \eps \tilde \beta})
\end{align}
>From formulas \eqref{formula beta}, \eqref{integral v f Z_a}  we see that
$$
\beta^2 \le C ( 1 + \sinh(\sqrt\eps \tilde \beta)  + \tilde \beta^2 ) .
$$
Using this inequality, standard elliptic estimates and the facts that $v\le 1$ and $f_t$ is uniformly bounded we obtain \eqref{bound T}. From here we deduce that $T$ is continuous and compact. Indeed, for the latter assertion, note that if $\mathcal B$ be a bounded set of $X \times [0,\infty) \times [0,1]$ then $T(\mathcal B)$ is a bounded set in $C^{1,\mu}([0,1]) \times \R$ for any $\mu \in (0,1)$ and taking $\mu>\alpha$ it follows that this set is precompact in $C^{1,\alpha}([0,1]) \times \R$.
\qed

\medskip


\begin{lem}
\label{lemma unique eps=0 t=0}
The operator $T_{0,0}$ has a unique fixed point $(u_0,\beta_0)$ in $X$.
\end{lem} 
\proof[\bf Proof.]
In this situation $u+v\equiv 1$ and hence the system reduces to
$$
-\Delta u = \beta^2 (1-u)(u-\theta)^+ \quad \hbox{in $B_1$}, \quad u|_{\partial B_1} = \theta
$$
with the additional requirement that
$$
- \pd{u}{\nu} = \theta.
$$
Let $w=u-\theta$. Then the equation for $w$ becomes
\begin{align}
\label{logistic}
-\Delta w = \beta^2 (1-\theta - w) w^+  \quad \hbox{in $B_1$}, \quad w|_{\partial B_1} = 0.
\end{align}
This equation is of logistic type and many properties are well known (see \cite{cantrell-cosner-book,hutson-lopezgomez-mischaikow-vickers}): 
\begin{enumerate}
\item
Any solution $w$ to \eqref{logistic} satisfies $0\le w < 1-\theta$, and either $w\equiv 0$ or $w>0$ in $B_1$.
\item
\eqref{logistic} has a nontrivial solution if and only if $\beta>\beta^*$ where $\beta^*$ is such that $(\beta^*) ^2 (1-\theta) = \lambda_1$ and  $\lambda_1$ denotes the first Dirichlet eigenvalue of $-\Delta$ in the unit ball.
\item
For $\beta>\beta^*$ there is a unique non-trivial solution, which we write as $w_\beta$.
\item
$w_\beta$ is monotone increasing with respect to $\beta$ and 
$$
\lim_{\beta\to\infty} w_\beta = 1-\theta
$$
uniformly on compact sets of $B_1$ and
$$
\lim_{\beta \to \beta^*} w_\beta = 0.
$$
\end{enumerate}
>From the above properties it follows that
$$
\lim_{ \beta\to\infty} - \pd{w_\beta}{\nu} = \infty \quad
\hbox{and} \quad
\lim_{ \beta\to0} - \pd{w_\beta}{\nu} = 0 .
$$
Moreover 
$$
\beta \mapsto - \pd{w_\beta}{\nu} \quad \hbox{is strictly increasing}.
$$
It follows that there is a unique $\beta>\beta^*$ such that 
$$
- \pd{w_\beta}{\nu} = \theta.
$$
We call this value $\beta_0$ and the let $u_0  = w_{\beta_0}$. Then $(u_0,\beta_0) \in X $ is the fixed point of $T_{0,0}$.
\qed


\begin{lem}
\label{sol is nondegenerate}
Let $(u_0,\beta_0)$ denote the unique fixed point of $T_{0,0}$ found in Lemma~\ref{lemma unique eps=0 t=0}. Then $(u_0,\beta_0)$ is nondegenerate.
\end{lem}
\proof[\bf Proof.] 
Let us write $(u_0,\beta_0) \in X$ the solution of $I-T_{0,0}=0$ of Lemma~\ref{lemma unique eps=0 t=0}. We have to verify that the linearization of $I-T_{0,0}$ around this solution is an invertible operator. 
Note that the operator $T_{0,0}$ involves only the nonlinearity $f_0(s) = (s-\theta)^+$
and that the fixed point $(u_0,\beta_0)$ from Lemma~\ref{lemma unique eps=0 t=0} satisfies $u_0>0$ in $[0,1)$ and $u_0'(1)<0$. 
Therefore for $\varphi$ in a neighborhood of $u_0$ in the topology of $C^{1,\alpha}([0,1])$ we have 
$f_0(u_0 + \varphi +\theta ) = (u_0 + \varphi )^+ = 
u_0 + \varphi$. 
Hence $T_{0,0}$ is differentiable at $(u_0,\beta_0)$.

We next compute the derivative of $T_{0,0}$ at $(u_0,\beta_0)$ 
in the direction of $(\varphi,\sigma) \in X$, which we write as 
$$
(\psi,\gamma ) = D T_{0,0}(u_0,\beta_0) (\varphi,\sigma).
$$
Then
\begin{align}
\label{psi gamma}
\psi = \psi_2+ \sigma \psi_1,
\quad
\gamma = \gamma_2 + \sigma \gamma_1
\end{align}
where $(\psi_2,\gamma_2)$ and $(\psi_1,\gamma_1)$ are computed as follows.

To compute $\psi_2$ we have linearize \eqref{T:sol u} with respect to $\tilde u$ and then set $\tilde u = u = u_0$, $\tilde \beta = \beta = \beta_0$. Since for $t=0$, we have $f_0(u_0 + \theta ) = u_0$ and $f_0'(u_0  + \theta ) = 1$,
the function $\psi_2$ satisfies
\begin{align}
\label{psi2}
\left\{
\begin{aligned}
&
-\Delta \psi_2 = 2 \beta_0 \gamma_2 v_0 u_0 
+ \beta_0^2 \pd{v}{\tilde u} u_0
+\beta_0^2 v_0 \varphi + 2 \beta_0 \gamma_2
\quad \text{in } B_1
\\
&
\psi_2(1) = 0, \quad \psi_2'(1)=0 ,
\end{aligned}
\right.
\end{align}
where $\gamma_2$ is adjusted so that $\psi_2$ satisfies both boundary conditions, $\pd{v}{\tilde u}$ is found solving
\begin{align}
\left\{
\label{dvdu}
\begin{aligned}
& 
\Delta \pd{v}{\tilde u} = 
\beta_0^2 \pd{v}{\tilde u} u_0 + \beta_0^2 v_0 \varphi
\quad\text{in } B_1
\\
&
\pd{v}{\tilde u}(1) + \pd{}{r}\pd{v}{\tilde u}(1) = 0,
\end{aligned}
\right.
\end{align}
and $v_0$ is the solution of \eqref{T:sol v} with $\tilde \beta = \beta_0$, $\tilde u = u_0$ and $t=0$. As explained in the proof of Lemma~\ref{lemma unique eps=0 t=0} then 
$$
v_0  = 1 - \theta - u_0.
$$
To compute $\psi_1$ we linearize \eqref{T:sol u} with respect to $\tilde \beta$ and then set $\tilde u = u = u_0$, $\tilde \beta = \beta = \beta_0$. Therefore $\psi_1$ satisfies
\begin{align}
\label{psi1}
\left\{
\begin{aligned}
&
- \Delta \psi_1 = 
2 \beta_0 \gamma_1 v_0 u_0
+ \beta_0^2 \pd{v}{\tilde \beta} u_0
+ 2 \beta_0 \gamma_1 - 2 \beta_2
\quad\text{in } B_1
\\
&
\psi_1(1) = 0 , \quad \psi_1'(1)=0,
\end{aligned}
\right.
\end{align}
where $\gamma_1$ is chosen so that both boundary conditions are satisfied and $\pd{v}{\tilde \beta}$ is computed from the equation
\begin{align}
\label{dvdb}
\left\{
\begin{aligned}
&
\Delta \pd{v}{\tilde \beta} 
= 2 \beta_0 v_0 u_0 + \beta_0^2 \pd{v}{\tilde \beta} u_0
\quad\text{in } B_1
\\
&
\pd{v}{\tilde \beta}(1) + \pd{}{r} \pd{v}{\tilde \beta}(1) = 0.
\end{aligned}
\right.
\end{align}

We need to find the kernel of 
$I - D_{\tilde u ,\tilde \beta}T_{0,0}(u_0,\beta_0)$, that is solutions $(\varphi,\sigma) \in X $ of:
\begin{align}
\label{eig problem}
D_{\tilde u,\tilde \beta} T_{0,0} (u_0,\beta_0) (\varphi,\sigma) = (\varphi,\sigma)
\end{align}
Using the notation above we may write $\varphi  = \psi_2 + \sigma \psi_1$, $\sigma = \gamma_2 + \sigma \gamma_1$. 
We claim that 
\begin{align}
\label{eq phi}
\varphi =  - \Big(
\pd{v}{\tilde u}+ \sigma
\pd{v}{\tilde \beta}
\Big) .
\end{align}
To prove this, let $z = \varphi + \pd{v}{\tilde u}+ \sigma
\pd{v}{\tilde \beta}$. Then a calculation shows that
$$
\Delta z = 0 \quad \text{in } B_1
$$
and $z(1) + z'(1) = 0$. This implies $z=0$ and 
\eqref{eq phi} follows.
>From the formulas \eqref{psi gamma}--\eqref{dvdb} we get
\begin{align*}
- \Delta \varphi  & = \beta_0^2 
\left(    1-\theta-2u_0 \right) \varphi   +  2 \sigma \beta_0 (1-\theta-u_0) u_0
\end{align*}
But also $\varphi$ satisfies
$$
\varphi(1) = 0,\quad
\varphi'(1)=0.
$$
The linear elliptic operator
$- \Delta - \beta_0^2  \left( 1-\theta-2 u_0\right) $ 
with Dirichlet boundary condition
arises as the linearization of \eqref{logistic} around the solution $u_0$, and it is well known 
that it has a positive first eigenvalue. Since $ \beta_0 (1-\theta-u_0) u_0  > 0$, for positive $\sigma$ we would have $\varphi>0$ and $\varphi'(1)<0$ and for negative $\sigma$ we would have $\varphi<0$ and $\varphi'(1)>0$, which are both impossible since $\varphi'(1)=0$. Hence $\sigma=0$ and also $\varphi=0$. This shows that the linear equation $D_{\tilde u,\tilde \beta} T_{0,0} (u_0,\beta_0) (\varphi,\sigma) =  (\varphi,\sigma)$ admits only the trivial solution. Since $D_{\tilde u,\tilde \beta} T_{0,0} (u_0,\beta_0)$ is compact we see that $I -  D_{\tilde u,\tilde \beta} T_{0,0} (u_0,\beta_0)$ is an isomorphism.

\qed

\proof[\bf Proof of Theorem~\ref{thm existence1}]

By Lemma~\ref{nonexistence small epsilon beta2} there exists $\delta>0$ and $\beta_1>0$ such that for all $\beta \ge \beta_1$, $\beta \sqrt \eps \le \delta$ and all $t \in [0,1]$ there is no solution in $\mathcal S$ to the system \eqref{system01}. This means that $T_{\eps,t}$ has no fixed point $(u,\beta)$ if $\beta \ge \beta_1$ and $\beta \sqrt \eps \le \delta$. Let 
$$
R_1 = C( 1+\sinh(\delta) + \beta_1^2)
$$
and define
\begin{align*}
\Omega_1 & = \{ (u,\beta) \in X: \| u \|_{C^{1,\alpha}([0,1])} < R_1 , \;   0 <\beta< \beta_1 \} .
\end{align*}
Let $\eps_1 = (\delta/\beta_1)^2$. Then $\Omega_1$ is a bounded open set of $X$ and for $0\le \eps\le \eps_1$ and $0\le t\le 1 $ the operator $T_{\eps,t}$ has no fixed point in $\partial \Omega_1$. Indeed, suppose $(u,\beta) \in \partial \Omega_1$ is a fixed point  of $T_{\eps,t}$. It is not possible that $\beta=\beta_1$ by Lemma~\ref{nonexistence small epsilon beta2}. The case $\beta = 0$ is also impossible. This means that $0<\beta<\beta_1$ and $\|u\|_{C^{1,\alpha}([0,1])} = R_1$. But by inequality \eqref{bound T} we would have
$$
R_1 = \|u\|_{C^{1,\alpha}([0,1])} \le  C( 1+\sinh(\sqrt\eps\beta) + \beta^2) < C( 1+\sinh(\delta) + \beta_1^2) = R_1
$$
which is impossible.
Thus  the Leray-Schauder degree $ deg(I - T_{\eps,t},\Omega_1,0) $ is well defined for $0\le \eps\le\eps_1$, $t \in [0,1]$. By Lemmas~\ref{lemma unique eps=0 t=0} and \ref{sol is nondegenerate} we have
$$
 deg(I - T_{0,0},\Omega_1,0)  = \pm 1
$$
and by homotopy invariance
$$
deg(I - T_{\eps,t},\Omega_1,0)  = \pm 1 \quad  \forall 0 \le \eps \le \eps_1 , \quad \forall t \in [0,1].
$$
This shows $T_{\eps,t}$ has at least one fixed point in $\Omega_1$ for $ 0 \le \eps \le \eps_1$ and $t\in [0,1]$. In particular the system \eqref{system01} has a solution with bounded $\beta$ for any $0 \le \eps \le \eps_1$.

We know there exists $\eps^*>0$ such that for $\eps>\eps^*$ the system \eqref{system01} has no solution, and hence $T_{\eps,1}$ has no fixed points in $X$ for such $\eps$. Let $0<\eps_0 < \eps_1$ and define
\begin{align*}
\Omega_2 & = \{ (u,\beta) \in X: \| u \|_{C^{1,\alpha}([0,1])} < R_2 , \;   0 <\beta< \beta_2 \}
\end{align*}
where
$$
R_2 = C( 1+\sinh( \sqrt {2 \eps^*} \beta_2) + \beta_2^2) \quad \hbox{and}\quad \beta_2 =  \theta /\sqrt\eps_0.
$$
Then $T_{\eps,t}$ has no fixed points on $\partial \Omega_2$ for $\eps_0 \le \eps \le 2 \eps^*$ and $t \in [0,1]$. Then using a homotopy along $\eps \in [\eps_0 , 2 \eps^*] $ we find
$$
deg(I - T_{\eps,1},\Omega_2,0)  = 0\quad \forall \eps_0 \le \eps \le 2 \eps^*.
$$
This implies that for $\eps_0 \le \eps \le \eps_1$ the system has at least another solution and that this solution has $\beta$ in the range $\delta/\sqrt \eps \le \beta \le \theta/\sqrt\eps$. Since $\eps_0$ is arbitrary we obtain the same conclusion for $0<\eps\le \eps_1$.
\qed

\section{Existence when $f$ is discontinuous}
\label{sect exsit discont}
In this section we prove Theorem~\ref{thm existence2}.
We assume that $f$ satisfies \eqref{h1}, \eqref{h1b}, \eqref{h2}
and that 
$$
\eta = \lim_{u\to \theta^+} f(u)>0.
$$
For $n\ge 1$ define
$$
f_n (t) = \min( f(t) , n(t-\theta)^+ ) .
$$
Then $f_n$ is continuous and satisfies \eqref{h1} and \eqref{h2} with a fixed constant. Moreover there is $h_0$ satisfying \eqref{f0 1} and \eqref{f0 2} and such that $f_n \ge h_0$ for all $n$. Such $h_0$ can be taken for instance as
$$
h_0(t) = \min( (t-\theta)^+ , \sigma)
$$
for some $\sigma>0$ small enough.

By Theorem~\ref{thm existence1} there is $\eps_1>0$ such that the system \eqref{system01} admits 2 solutions $(u_n^1,v_n^1)$,   $(u_n^2,v_n^2)$ for all $\eps \in (0,\eps_1)$. Moreover there is a fixed number $\beta_1$ such that  $\beta_n^1 \le \beta_1$ and $ \beta_n^2 \in (\delta/\sqrt\eps,1/(\theta \sqrt \delta)) $ for all $n\ge 1$ and $\eps \in (0,\eps_1)$. From now on we fix $\eps \in (0,\eps_1)$ and study the limit as $n\to \infty$ of any of the 2 solutions which we call just $(u_n,v_n)$ with parameter $\beta_n$. Since $\beta_n$ is bounded and by standard elliptic estimates we may assume that $\beta_n \to \beta$, $u_n \to u$ and $v_n \to v$ in $C^{1,\alpha}([0,1])$.

Since $f_n(u_n)$ is uniformly bounded we may extract a further subsequence such that $f_n(u_n) \to \varphi $ weakly-* in $L^\infty(0,1)$. The function $\varphi$ then satisfies $0 \le \varphi \le \max_{[0,1]} f$ in $[0,1]$. Since $u_n$ satisfies a linear equation for $r>1$ and vanishes at infinity it has an explicit formula
$$
u_n(r) = \theta \frac{e^{-(r-1)\beta_n \sqrt \eps}}{r} \quad \hbox{for all $r\ge 1$.}
$$
Setting $\varphi(r) = 0$ for $r>1$ we see that  $u,v$ satisfy
\begin{align*}
\left\{
\begin{aligned}
\Delta u - \eps \beta^2 u &= - \beta^2 v \varphi \quad \hbox{in $\R^3$}
\\
\Delta v &=  \beta^2 v \varphi \quad \hbox{in $\R^3$}
\\
\lim_{|x|\to\infty} u(x) & = 0  , \qquad 
\lim_{|x|\to\infty} v(x) = 1 
\end{aligned}
\right.
\end{align*}
The set $D =\{ r \in [0,\infty) \ : \ u(r) \not=\theta \}$ is open (relative to $[0,\infty))$ and for any $r \in D$ we have $f_n(u_n(r)) \to f(u(r))$ as $n\to\infty$. It follows that $\varphi = f(u)$ in $D$.  Outside $D$ we have $\varphi \le \eta$ a.e.
Indeed, let $\psi \in C_0^\infty((0,1))$, $\psi\ge 0$. Then
$$
\int_0^1 f_n(u_n) \psi = \int_{D \cap [0,1]} f_n(u_n) \psi +  \int_{D^c \cap [0,1]}  f_n(u_n) \psi.
$$
Since $f_n(u_n) \to f(u)$ in $D$ by dominated convergence we have
$$
\int_{D \cap [0,1]} f_n(u_n) \psi  \to \int_{D \cap [0,1]} f(u) \psi .
$$
On the other hand, if $r\in D^c$ $\limsup_{n\to+\infty} f_n(u_n(r)) \le \eta$ so that $ (f_n(u_n) - \eta)^+ \to 0$ on $D^c$. Then
$$
\int_{D^c \cap[0,1]} (f_n(u_n) -\eta) \psi = \int_{D^c \cap[0,1]}  (f_n(u_n) -\eta)^+  \psi - \int_{D^c \cap[0,1]}  (f_n(u_n) -\eta)^-  \psi \le o(1)
$$
where $o(1)$ denotes a sequence converging to 0 as $n\to+\infty$. It follows that 
$$
\int_{D^c \cap [0,1]}  f_n(u_n) \psi \le \int_{D^c \cap [0,1]} \eta \psi + o(1)
$$
and hence
$$
\int_0^1 \varphi \psi = \lim_{n\to +\infty} \int_0^1 f_n(u_n) \le \int_{D\cap [0,1]} f(u) \psi + \int_{D^c \cap [0,1]} \eta \psi.
$$
This shows that $\varphi \le \eta$ a.e.\@ in $D^c$.

Our main objective is to show that the complement of $D$ is finite.
If $u'(r) \not = 0$ whenever $u(r) = \theta$ then $D^c$ is discrete and since it is contained in $[0,1]$ it is finite. 

Let us analyze the case where for some $r \in [0,1] $ we have $u(r) = \theta$ and $u'(r)=0$.
Let 
$$
r_0 = \sup \ \{ r>0 \ :  \ u(r) = \theta, \ u'(r) =0  \}.
$$
Then $r_0<1$ and $u(r_0) = \theta$, $u'(r_0) = 0$. We assert that there is a small interval $(r_0,r_0+\sigma)$, $\sigma>0$ such that $u>\theta$ in $(r_0,r_0+\sigma)$. To prove this we start ruling out the possibility that $u(r_n) = \theta$ for some infinite sequence $r_n \searrow r_0$. We actually may assume that if $n$ is even then  $u>\theta$  on  $(r_{n+1},r_n)$ and if $n$ is odd  $u<\theta$  on  $(r_{n+1},r_n)$. Let us see that the following holds
\begin{align}
\label{derivatives}
u'(r_{n+1}) = - u'(r_n) + O(\ell_n^2)
\end{align} 
where
$$
\ell_n = r_{n+1} - r_n
$$
and  $O(\ell_n^2)$ denotes a sequence bounded by $C \ell_n^2$ with $C$ independent of $n$ as $n\to \infty$.
Suppose first that $u<\theta$  on  $(r_{n+1},r_n)$ and define $\tilde u_n$ by
$$
u(r) = \theta + \ell_n^2\tilde u_n(( r - r_{n+1})/\ell_n) \quad r \in [r_{n+1} ,r_n].
$$
Then
$$
\tilde u_n'' + \frac{2\ell_n}{r \ell_n + r_{n+1} } \tilde u_n' = \eps\beta^2 ( \ell_n^2\tilde u_n + \theta) \quad \hbox{for } r \in [0,1].
$$
$$
\tilde u_n(0) = \tilde u_n(1) = 0.
$$
This equation implies that $\tilde u_n$ and $\tilde u_n'$ are uniformly bounded on $[0,1]$. Therefore
$$
\tilde u_n'' = \eps\beta^2  \theta+ O( \ell_n)  \quad \hbox{for } r \in [0,1] 
$$
hence integrating
$$
\tilde u_n(r) = \frac{1}{2} \eps\beta^2 ( 1 + O(\ell_n)) r( r-1) \quad \hbox{for } r \in [0,1] 
$$
and
\begin{align}
\label{der rn}
\tilde u_n'(r) = \frac{1}{2} \eps\beta^2 ( 1 + O(\ell_n)) ( 2r-1) \quad \hbox{for } r \in [0,1] .
\end{align}
This proves \eqref{derivatives} in this case.

Now suppose that $u>\theta$  on  $(r_{n+1},r_n)$. Defining $\tilde u_n$ as before we now obtain the following equation
$$
\tilde u_n'' + \frac{2\ell_n}{r \ell_n + r_{n+1} } \tilde u_n' = \eps\beta^2 ( \ell_n^2\tilde u_n + \theta)  - \beta^2 \tilde v_n f(\theta + \ell_n^2 \tilde u_n) \quad \hbox{for } r \in [0,1].
$$
$$
\tilde u_n(0) = \tilde u_n(1) = 0
$$
where
$$
\tilde v_n(r) = v ( r_{n+1} + \ell_n r) \quad r \in [0,1].
$$
This equation implies that $\tilde u_n$, $\tilde u_n'$ are uniformly bounded on $[0,1]$.
Thus
$$
\tilde u_n''= \eps\beta^2  \theta - \beta^2 \tilde v_n(0) f(\theta + \ell_n^2 \tilde u_n) + O(\ell_n) \quad \hbox{for } r \in [0,1].
$$
Multiplying this equation by $\tilde u_n$ and integrating in $[0,1]$ yields
$$
\frac{1}{2} (\tilde u_n(1)')^2 -  \frac{1}{2} (\tilde u_n(0)')^2 = O(\ell_n)
$$
which shows that 
$$
\tilde u_n(1)' = - \tilde u_n(0)' + O(\ell_n)
$$
and proves \eqref{derivatives} in this case.

Using \eqref{derivatives} inductively we find
$$
u'(r_n) = (-1)^k u'(r_{n+k}) + \sum_{j=n}^{n+k-1} O(\ell_j^2) \quad \forall n \ge 1, \  k \ge 1.
$$
Letting $k\to+\infty$, and using that $u'(r_0)=0$ yields
$$
u'(r_n) =  \sum_{j=n}^{\infty} O(\ell_j^2) \quad \forall n \ge 1.
$$
Choose a subsequence $n_i \to \infty$ as $i\to\infty$ such that $\ell_{n_i} \ge \ell_{j}$ for all $j \ge n_i$. Then
$$
|u'(r_{n_i}) | \le C \ell_{n_i} \sum_{j=n_i}^\infty \ell_j.
$$
But \eqref{der rn} implies that $|u'(r_{n}) | = \ell_n \eps\beta^2 (1+O(\ell_n))/2$. It follows that 
$$
\ell_{n_i} \eps\beta^2 (1+O(\ell_{n_i} ))/2 \le  C \ell_{n_i} \sum_{j=n_i}^\infty \ell_j
$$
which is impossible as $i\to\infty$.

This establishes the assertion that $u$ cannot oscillate infinitely many times around $\theta$ to the right of $r_0$. It follows that for some $\sigma>0$ either $u>\theta$ in $(r_0,r_0+\sigma)$ or $u<\theta$ in $(r_0,r_0+\sigma)$. The latter can in fact not occur by the Hopf lemma.

Since $u$ satisfies the ODE $u''+\frac{2}{r} u' = \eps\beta^2 u - \beta^2 v f(u)$ in $(r_0,r_0+\sigma)$ and in this interval $u>\theta$ we see that $u$ is $C^2 $ in $[r_0,r_0+\sigma)$. Since $r_0$ is a minimum of $u$ restricted to $[r_0,r_0+\sigma)$ it follows that $\lim_{r\to r_0^+} u''\ge0$. Which yields the following inequality
\begin{align}
\label{ineq r0}
\eps \beta^2 \theta  - \beta^2 v(r_0) \eta \ge 0.
\end{align}
Now, if $v$ is constant in $[0,r_0]$ this means that $\varphi = 0$ a.e.\@ in $[0,r_0]$ and $u \le \theta$ in $[0,r_0]$. Actually $u<\theta$ in $[0,r_0)$ by the strong maximum principle and this finishes the proof in this case. If $v$ is not constant, then $v(r)<v(r_0)$ for all $r \in [0,r_0)$. Let us verify that $D^c$ is discrete. Suppose that for some $r_1 \in [0,r_0)$ we have $u(r_1) = \theta$. If $u'(r_1)\not=0$ then $r_1$ is isolated. If $u'(r_1)=0$ then integrating the equation for $u$ we have
$$
r^2 u'(r) = \int_{r_1}^r t^2 ( \eps \beta^2 u - \beta^2 v \varphi) \, d t.
$$
For $r$ near $r_1$ we have
\begin{align*}
& \hbox{if } u(r)<\theta \hbox{ then } \eps \beta^2 u(r) - \beta^2 v(r) \varphi(r) > 0 
\\
& \hbox{if } u(r) = \theta \hbox{ then } \eps \beta^2 u(r) - \beta^2 v(r) \varphi(r) >  \eps \beta^2 \theta - \beta^2 v(r_0) \eta \ge 0
\\
& \hbox{if } u(r)>\theta \hbox{ then } \eps \beta^2 u(r) - \beta^2 v(r) \varphi(r)  = \eps \beta^2 u(r) - \beta^2 v(r) f(u(r)) .
\end{align*}
In the second line above we have used that $v(r)<v(r_0)$,  $\varphi \le \eta $ and \eqref{ineq r0}. In the third line above we may say that $\eps \beta^2 u(r) - \beta^2 v(r) f(u(r)) >0$ if $r$ is sufficiently close to $r_1$ by continuity, $v(r)<v(r_0)$ and \eqref{ineq r0} (we may regard $f$ as continuous here since we are working with values above $\theta$). This shows that $u$ is strictly convex in a neighborhood of $r_1$ and hence there are no other points in $D^c$ close to $r_1$. This shows that $D^c$ is discrete, hence finite, and finishes the proof of the theorem.

%



\section{The Heaviside ignition}
\label{The Heaviside ignition}

In this section we perform explicit calculations for the ignition nonlinearity $f(u) = \chi_{[u>\theta]}$, where $0<\theta<1$. We first reduce the differential equations to a finite number of equations in some parameters. In a second part we obtain rigorously the bifurcation diagram of Figure~\ref{fig1} for small $\eps$. 

We rewrite \eqref{system01} in the form
\begin{align}
\label{p3}
\left \{
   \begin{array}{l}
       \Delta u -\beta^2\eps u = 0 \quad \hbox{in  $(1,+\infty)$}\\
       \Delta v =0  \quad \hbox{in  $(1,+\infty)$}\\
       u(1)=\theta,\ \ \ v(1)=\gamma \\
       u(+\infty)= 0 \ \ \ v(+\infty)=1 \\
   \end{array}
  \right.
\end{align}
and
\begin{align}
 \left \{
   \begin{array}{l}
        \Delta u -\beta^2\eps u=-\beta^2vf(u)  \quad \hbox{in  $(0,1)$} \\
        \Delta v =\beta^2vf(u)  \quad \hbox{in  $(0,1)$} \\
       u(1)=\theta \ \ \ v(1)=\gamma \\
           u^{\prime}(0)=0  \ \ \ v^{\prime}(0)=0 \\
u^{\prime}(1^-)=u^{\prime}(1^+) \ \ \  v^{\prime}(1^-)=v^{\prime}(1^+)           
   \end{array}
  \right.\label{p4}
\end{align}
where $\gamma$ is a parameter to be adjusted.

Observe that problem \eqref{p3} can explicitly be solved. Indeed, the following functions solve \eqref{p3}
\begin{align*}
u=\frac{\theta e^{-\beta\sqrt{\eps}(r-1)}}{r} 
\qquad v=1-\frac{1-\gamma}{r}.
\end{align*}
Hence the last condition in \eqref{p4} becomes
\begin{align*}
u^{\prime}(1^+)=-\theta(1+\beta\sqrt{\eps})
\qquad
v^{\prime}(1^+)=1-\gamma.
\end{align*}

\subsection{Solving \eqref{p4} with the assumption $u>\theta$ in $B(0,1)$}
\label{sub reaction ball}
Assuming that $u>\theta$, in $B(0,1)$ the subsystem \eqref{p4} is reduced to the following linear problem
 \begin{align}
        &\Delta u -\beta^2\eps u=-\beta^2v\qquad \text{ in } B(0,1) \label{p5-l1}\\
        &\Delta v =\beta^2v \qquad \text{ in } B(0,1) \label{p5-l2}\\
       &u(1)=\theta \ \ \ v(1)=\gamma \label{p5-l3}\\
    &v^{\prime}(1^-)=1-\gamma.\label{p5-l4}\\
      & u^{\prime}(1^-)=-\theta(1+\beta\sqrt{\eps}) \label{p5-l5}
\end{align}
For $\eps \neq 1$, observe that the two following functions
\begin{align}
\label{v small beta inner}
&v:=\frac{1}{\beta ch(\beta)} \frac{sh(\beta r)}{r} \quad r\in(0,1) 
\\
\label{u small beta inner}
& 
u:= - \frac{v}{1-\eps} +\frac{\left(\theta + \frac{th(\beta)}{\beta(1-\eps)}\right)}{sh(\beta\sqrt{\eps})}\frac{sh(\beta\sqrt{\eps}r)}{r} \quad  r\in(0,1)
\end{align}
solve equations \eqref{p5-l1}-\eqref{p5-l4}.
>From the formula for $v$ we can see that 
$$
\gamma=\frac{th(\beta)}{\beta}.
$$
Now to obtain a solution of subsystem \eqref{p4}, it then remains to adjust $\beta$ and $\eps$ to get \eqref{p5-l5}.
This equation is
$$ 
-\theta(1+\beta\sqrt{\eps})= u^{\prime}(1^-)=\frac{v'(1)}{1-\eps} +\frac{\left(\theta -\frac{th(\beta)}{\beta(1-\eps)}\right)}{sh(\beta\sqrt{\eps})}(\beta\sqrt{\eps} ch(\beta\sqrt{\eps})-sh(\beta\sqrt{\eps}))
$$
which after simplification yields
\begin{equation}\label{implicitsol}
\theta(1-\eps)(1+th(\beta\sqrt{\eps}))+\frac{th(\beta)}{\beta} = \frac{th(\beta\sqrt{\eps})}{\beta\sqrt{\eps}} .
\end{equation}
Let us define 
\begin{align}
\label{function g}
g(\eps,\beta):= \theta(1-\eps)(1+th(\beta\sqrt{\eps}))+\frac{th(\beta)}{\beta} -\frac{th(\beta\sqrt{\eps})}{\beta\sqrt{\eps}}.
\end{align}
It can be shown that there is $\eps^*>0$ such that for $\eps \in (0,\eps^*)$, $g(\eps,\cdot)$ has 2 zeros, $\beta_\eps^-<\beta_\eps^+$. 
The function $g(0,\cdot)$ has a unique zero $\beta_0 = \lim_{\eps\to0} \beta_\eps^-$. 
A calculation (see Lemma~\ref{u>theta small eps}) shows that for $\eps>0$ small and $\beta = \beta_\eps^-$ the solution $u$ obtained by the formula \eqref{u small beta inner} satisfies $u>\theta$ in $[0,1)$. Along the lower branch $\beta_\eps^-$, $0\le \eps \le \eps^*$ we observe numerically that $u>\theta$ in $[0,1)$. After reaching $\eps^*$ we continue the curve along the upper branch $\beta_\eps^+$ now decreasing $\eps$ from $\eps^*$ to 0 until we reach a critical value $\eps=\eps_0$ where the corresponding solution given by  \eqref{u small beta inner} satisfies $u>\theta$ in $(0,1)$ but $u(0) = \theta$.

\subsection{Solving the subsystem \eqref{p4} with the assumption $u>\theta$ in $B(0,1)\setminus B(0,\eta)$ and $u\le \theta$ in $B(0,\eta)$ for some $\eta>0$} \
\label{analysis reaction annulus}

Under this assumption the problem \eqref{p4} decouples again into two subsystems
\begin{align}
\left \{
   \begin{array}{ll}
        \Delta u -\beta^2\eps u=-\beta^2v & \text{in } B(0,1)\setminus B(0,\eta) \\
        \Delta v =\beta^2v & \text{in } B(0,1)\setminus B(0,\eta) \\
       u(1)=\theta=u(\eta) & v(1)=\gamma \ \ \ v(\eta)=\delta \\
       u^{\prime}(1^-)=u^{\prime}(1^+) &  v^{\prime}(1^-)=v^{\prime}(1^+)\\
       u^{\prime}(\eta^-)=u^{\prime}(\eta^+) &  v^{\prime}(\eta^-)=v^{\prime}(\eta^+)\\
   \end{array}
  \right.\label{p7}
\end{align}
and
\begin{align}
  \left \{
   \begin{array}{ll}
        \Delta u -\beta^2\eps u=0& \text{ in }  B(0,\eta) \\
        \Delta v =0& \text{ in }  B(0,\eta) \\
       u(\eta)=\theta \ \ \ v(\eta)=\delta
   \end{array}
  \right.\label{p8}
\end{align}
where $\delta,\gamma,\eta$ are parameters to be found. Using the conditions at $r=\eta$ and $r=1$ we will be able to reduce the parameters to only $\beta$ and $\eta$ which will be implicitly defined as functions of $\eps$ by 2 equations.


Observe that since $v$  is harmonic  in $B(0,\eta)$ with the Dirichlet condition  $v(\eta)=\delta$, we deduce that $v\equiv \delta$ and therefore $v^{\prime}(\eta-)=0$.

Let us now  solve \eqref{p7} explicitly.  First, observe that the function $v=A\, \frac{sh(\beta r)}{r}+B\, \frac{ch(\beta r)}{r}$ satisfies
$$\Delta v-\beta^2 v=0 \quad \text{ in }\quad B(0,1)\setminus B(0,\eta).$$
Therefore, choosing $A$ and $B$ such that
\begin{align*}
&\gamma = A\, sh(\beta)+ B\, ch(\beta)\\
&1-\gamma = A[\beta ch(\beta) -sh(\beta)] +B[\beta sh(\beta) -ch(\beta)],
 \end{align*}
it follows that $v$ satisfies the right boundary conditions on $\partial B(0,1)$.
Solving for $A$ and $B$ we find
\begin{align*}
A &= \frac{ch(\beta)-\gamma\beta sh(\beta)}{\beta}
\quad\hbox{and}\quad
B =  \frac{sh(\beta)-\gamma\beta ch(\beta)}{\beta}.
\end{align*}
To obtain a full solution it remains to adjust $\gamma$ in such way that $v'(\eta)=0$.
But
$$
v^{\prime}=A[\frac{\beta ch(\beta r)}{r} -\frac{sh(\beta r)}{r^2}] +B[\frac{\beta sh(\beta r)}{r} -\frac{ch(\beta r)}{r^2}] 
$$
and setting  $v^{\prime}(\eta)=0$ in the above equation  yields
\begin{align*}
 &A\left[\frac{\beta ch(\beta \eta)}{\eta} -\frac{sh(\beta \eta)}{\eta^2}\right] +B\left[\frac{\beta sh(\beta \eta)}{\eta} -\frac{ch(\beta \eta)}{\eta^2}\right] =0 \\
&A\Big(\eta\beta ch(\beta \eta)-sh(\beta \eta)\Big) + B \Big( \eta\beta sh(\beta \eta)-ch(\beta \eta)\Big)=0
\end{align*}
Substituting now the values of $A$ and $B$ in the above equation, it follows
$$
\Big(ch(\beta)-\gamma \beta sh(\beta)  \Big)C - \Big(sh(\beta)-\gamma \beta ch(\beta)  \Big)D =0,
$$
where
\begin{align*}
C:=\Big(\eta\beta ch(\beta \eta)-sh(\beta \eta)\Big)\\
D:=\Big( \eta\beta sh(\beta \eta)-ch(\beta \eta)\Big)
\end{align*}
Therefore, we have
\begin{equation}
\gamma =\frac{C-Dth(\beta)}{\beta\Big( Cth(\beta)-D\Big)}\label{fbgamma}
\end{equation}
which after simplification reads
$$
\gamma
=\gamma(\beta,\eta)=
\frac{\eta\beta -th(\eta\beta)-th(\beta)\Big[ \eta\beta th(\beta\eta) -1\Big]}{\beta\left[\Big(\eta\beta- th(\eta \beta)\Big)th(\beta) +1 - \eta\beta th(\eta\beta) \right]} .
$$
Thus we may compute now $\delta = v(\eta)$:
$$
\delta=\left(\frac{ch(\beta)-\gamma\beta sh(\beta)}{\beta}\right) \frac{sh(\beta \eta)}{\eta} -\left( \frac{sh(\beta)-\gamma\beta ch(\beta)}{\beta}\right) \frac{ch(\beta \eta)}{\eta}.
$$
>From \eqref{p8} we know that  in $B(0,\eta)$, we have the following
$$
u= \frac{\theta\eta}{sh(\eta\beta\sqrt{\eps})}\frac{sh(r\beta\sqrt{\eps})}{r}.
$$
Therefore, 
\begin{align}
\label{u prime eta 0}
u^{\prime}(\eta)=\theta\left[\frac{\beta\sqrt{\eps}ch(\eta\beta\sqrt{\eps})}{sh(\eta\beta\sqrt{\eps})} -\frac{1}{\eta}\right]=\frac{\theta ch(\eta\beta\sqrt{\eps})}{\eta sh(\eta\beta\sqrt{\eps})}[\eta\beta\sqrt{\eps} -th(\eta\beta\sqrt{\eps})].
\end{align}

Let us now find the solution $u$ in $B(0,1)\setminus B(0,\eta)$.  As for the construction of $v$, let us observe that for $\eps\neq 1$,
\begin{align}
\label{formula u}
u:=\frac{1}{\eps - 1}v+E\frac{sh(r\beta\sqrt{\eps})}{r}+ F\frac{ch(r\beta\sqrt{\eps})}{r}
\end{align}
solves the following equation
$$
\Delta u -\beta^2\eps u = -\beta^2 v 
$$
By taking
\begin{align}
\label{E}
E(\beta,\eta,\gamma):=\frac{ch(\beta\sqrt{\eps})-\gamma\beta\sqrt{\eps}sh(\beta\sqrt{\eps})}{(1-\eps)\beta\sqrt{\eps}} - \theta\Big(ch(\beta\sqrt{\eps})+sh(\beta\sqrt{\eps}) \Big)\\
\label{F} 
F(\beta,\eta,\gamma):= - \frac{sh(\beta\sqrt{\eps})-\gamma\beta\sqrt{\eps}ch(\beta\sqrt{\eps})}{(1-\eps)\beta\sqrt{\eps}} +  \theta\Big(ch(\beta\sqrt{\eps})+sh(\beta\sqrt{\eps}) \Big)
\end{align}
we easily verify that  $u(1)=\theta$ and $u^{\prime}(1^-)=-\theta(1+\sqrt{\eps}\beta)$. It remains now to impose the boundary condition on $\partial B(0,\eta)$. From the formula for $u$ \eqref{formula u} we have
\begin{align}
\nonumber
u(\eta) &= \frac{\delta}{\eps - 1}+E\frac{sh(\eta\beta\sqrt{\eps})}{\eta}+ F\frac{ch(\eta\beta\sqrt{\eps})}{\eta}
\\
\label{u prime eta 01}
u'(\eta) & =E\beta\sqrt{\eps}\frac{ch(\beta \sqrt{\eps} \eta)} {\eta}+F\beta\sqrt{\eps}\frac{sh(\beta \sqrt{\eps}\eta)}{\eta}-E\frac{sh(\beta\sqrt{\eps}\eta)}{\eta^2}
 -F\frac{ch(\beta\sqrt{\eps}\eta)}{\eta^2}
\end{align}
since $v'(\eta)=0$. Thus the equation $u(\eta)=\theta$ becomes
\begin{align}
\label{u eta}
\frac{\delta}{\eps - 1}+E(\beta,\eta,\gamma(\beta,\eta))
\frac{sh(\eta\beta\sqrt{\eps})}{\eta}+ F(\beta,\eta,\gamma(\beta,\eta))\frac{ch(\eta\beta\sqrt{\eps})}{\eta} =\theta
\end{align}
and \eqref{u prime eta 0} combined with \eqref{u prime eta 01} becomes
\begin{align}
\nonumber
E(\beta,\eta,\gamma(\beta,\eta))\eta\beta\sqrt{\eps}ch(\beta \sqrt{\eps} \eta)+
F(\beta,\eta,\gamma(\beta,\eta))\eta\beta\sqrt{\eps}sh(\beta\sqrt{\eps}\eta)
\\
\nonumber
-E(\beta,\eta,\gamma(\beta,\eta)) sh(\beta\sqrt{\eps}\eta)-F(\beta,\eta,\gamma(\beta,\eta)) ch(\beta \sqrt{\eps} \eta) &
\\
\label{u prime eta}
=\theta\left[\beta\sqrt{\eps} coth(\eta\beta\sqrt{\eps}) -\frac{1}{\eta}\right] . &
\end{align}

We solve numerically 
equations \eqref{u eta} and \eqref{u prime eta}  with $\beta$ and $\eta$ as unknowns that depend on $\eps$, for $0<\eps<\eps_0$, where $\eps_0$ is the critical value of $\eps$ described in Subsection~\ref{sub reaction ball}.
The result from this numerical computation is shown in Figure~\ref{fig1}.

\subsection{Results for small $\eps$}
\begin{lem}
\label{u>theta small eps}
For $\eps>0$ suitably small $\beta_\eps^-$ is a smooth function of $\eps$ and the solution given by \eqref{u small beta inner} corresponding to $\beta = \beta_\eps^-$ satisfies $u>\theta$ in $[0,1)$.
\end{lem}
\proof[\bf Proof]
Let $\beta_1$ the unique solution of the equation
$$\frac{th(\beta)}{\beta}=\frac{(1-\theta)}{2}.$$
Then 
\begin{align*}
 g(\eps, \beta_1,\theta)&=\theta(1-\eps)(1+th(\beta_1\sqrt{\eps}))+\frac{1-\theta}{2} -\frac{th(\beta_1\sqrt{\eps})}{\beta_1\sqrt{\eps}}\\
&=\theta(1-\eps)th(\beta_1\sqrt{\eps}))+\frac{1+\theta -2\theta\eps}{2} -\frac{th(\beta_1\sqrt{\eps})}{\beta_1\sqrt{\eps}}\\
&=\theta(1-\eps)th(\beta_1\sqrt{\eps}))-\frac{1-\theta}{2}-\theta\eps +[1 -\frac{th(\beta_1\sqrt{\eps})}{\beta_1\sqrt{\eps}}]
\end{align*}
For  $\beta_1\sqrt{\eps}<1$,  using that $\beta_1\sqrt{\eps}-\frac{(\beta_1\sqrt{\eps})^3}{3}\le th(\beta_1\sqrt{\eps})\le \beta_1\sqrt{\eps}$, we end up with

\begin{align*}
g(\eps, \beta_1,\theta)&\le \theta(1-\eps)\beta_1\sqrt{\eps}-\frac{1-\theta}{2}-\theta\eps +\frac{(\beta_1\sqrt{\eps})^2}{3}\\
                                   &\le\frac{(\beta_1\sqrt{\eps})^2}{3}+ \theta\beta_1\sqrt{\eps}-\frac{1-\theta}{2}-\theta\eps.
\end{align*}
Computing now the positive roots of $    \frac{X^2}{3}+\theta X -\frac{1-\theta}{2}$, yields $g(\eps, \beta_1,\theta)<0$ for $\sqrt{\eps}\beta_1\le\frac{3\theta\left(\sqrt{ 1+\frac{2(1-\theta)}{3\theta^2}  }       -1\right)}{2}$.

Therefore,  we achieve  $g(\eps,\beta_1,\theta)< 0$, for $\eps$ in $ (0,\min\{\left(\frac{3\theta\left(\sqrt{ 1+\frac{2(1-\theta)}{3\theta^2}  }       -1\right)}{2\beta_1}\right)^2,\frac{1}{\beta^2_1}\})$.
Let us denote $$\eps_1:=\left(\frac{3\theta\left(\sqrt{ 1+\frac{2(1-\theta)}{3\theta^2}  }       -1\right)}{2\beta_1}\right)^2.$$
Hence, for any $(\eps,\theta)$ in $ (0,\min\{\eps_1,\frac{1}{\beta^2_1}\})\times(0,1)$ fixed there exists two possible solution to $g(\eps,\beta,\theta)=0$.
Moreover, $\beta^-<\beta_1<\beta^+$.

To obtain a  solution, to problem \eqref{p4}, we still need to show that the constructed solutions effectively satisfies the conditions $u>\theta$ in $B(0,1)$.
Let us observe that $u-\theta$ satisfies:
$$\Delta(u-\theta)-\eps\beta^2(u-\theta)=-\beta^2(v-\eps \theta).$$
If $v(0)-\eps \theta>0$ then using the maximum principle, it follows that $(u-\theta)>0$.
Since $v(r)=\left(\frac{1}{ch(\beta)\beta}\right) \frac{sh(\beta r)}{r}$, it follows that $v(0)=\frac{2}{ch(\beta)}$.
Hence, we end up with the condition
$$ \frac{2}{ch(\beta)}\ge \eps \theta.$$
We can conclude with a final estimate on $\eps$, namely since $\beta^-<\beta_1$, and $\frac{1}{ch(\beta)}$ is a decreasing function, we have
the following uniform estimates on $\eps$
$$\eps\le \frac{(1-\theta)\beta_1}{\theta sh(\beta_1)}.$$
We therefore have construct a solution to the problem (1), when $(\eps,\theta)$ belongs to $(0,\eps_0(\theta))\times(0,1)$ where
$$
\eps_0(\theta):=\min\{\frac{1}{\beta_1^2}; \frac{(1-\theta)\beta_1}{\theta sh(\beta_1)};\eps_1;\frac{(1-\theta)^2}{2}      \}
$$
\qed

\begin{lem}
For $\eps>0$ small the system \eqref{u eta}, \eqref{u prime eta} has a unique solution $\beta_\eps^+,\eta$ with $\eta>0$,  $\beta_\eps^+$ is a smooth function of $\eps$ and $u$ defined by \eqref{formula u} satisfies $u>\theta$ in $(\eta,1)$.
\end{lem}

\proof[\bf Proof]
The 2 equations \eqref{u eta} and \eqref{u prime eta} can be written in the form
\begin{align*}
E \sinh(a\eta) + F \cosh(a\eta) &= \eta \theta + \frac{\eta \delta}{1-\eps}
\\
E  \left( a  \eta \cosh(a \eta) - \sinh(a \eta) \right) + F \left( a \eta \sinh(a \eta) - \cosh(a \eta)  \right) &=  \eta \theta\left( \frac{a \eta}{\tanh(a\eta)} - 1 \right) 
\end{align*}
where
$$
a = \beta \sqrt \eps.
$$
Instead of $\eps$, $\beta$, $\eta$ consider the variables $t$, $a$, $x$ defined by the following relations
$$
t = \sqrt \eps, \quad a = \beta \sqrt \eps, \quad \eta = 1 - \frac{x}{\beta}.
$$
Consider also the function $G(t,a,x) = ( G_1(t,a,x), G_2(t,a,x) )$  where
\begin{align*}
G_1(t,a,x) &= \beta\left[ E \sinh(a\eta) + F \cosh(a\eta) - \eta \theta - \frac{\eta \delta}{1-\eps} \right] 
\\
G_2(t,a,x) &= E  \left( a  \eta \cosh(a \eta) - \sinh(a \eta) \right) + F \left( a \eta \sinh(a \eta) - \cosh(a \eta)  \right) 
\\
& \quad 
- \eta \theta\left( \frac{a \eta}{\tanh(a\eta)} - 1 \right) 
\end{align*}
where $\delta,\gamma,E,F$ are as before.
We see that we have a solution to the system if and only if $G(t,x,a)=0$.

Using the fact that
$$
\lim_{t \to 0} \beta \gamma =  \frac{1}{\tanh(x)} \quad \hbox{and} \quad \lim_{t \to 0} \beta \delta = \frac{1}{\sinh(x)}
$$
it is not difficult to verify that $G$ maybe extended in a $C^1$ manner for $t=0$, $a,x \in(0,\infty)$ with the values
\begin{align*}
G_1(0,a,x) &= - x + \theta a x + \frac{1}{\tanh(x)} + x \theta - \frac{1}{\sinh(x)} 
\\
G_2(0,a,x) &=  1 - \theta(1+a) - \theta\left( \frac{a }{\tanh(a)} - 1 \right) 
\end{align*}
Note that the equation  $G_2(0,a,x) =0$ is equivalent to
$$
\frac{\frac{a }{\tanh(a)} - 1}{ 1 - \theta(1+a)} = \frac{1}{\theta} 
$$
which is seen to have a unique solution $a_0 \in (0,1/\theta-1)$ since the left hand side defines a strictly increasing function of $a$ which goes to zero as $a \to 0$ and goes to $+\infty$ as $a \nearrow 1/\theta-1$.

Now the equation $G_1(0,a_0,x) = 0$ can be written as 
\begin{align}
\label{eq x}
1 - \frac{1}{x \tanh(x)} + \frac{1}{x \sinh(x)} = \theta (1 + a_0) \in (0,1).
\end{align}
The above left hand side is strictly increasing with limits $1/2$ as $x\to 0$ and 1 as $x \to \infty$. Using the inequality
$$
\frac{a}{\tanh(a)} - 1 \le a \quad \forall a\ge 0
$$
we see that $ 1- \theta - \theta a_0 \le \theta a_0$ which yields $a_0 \ge \frac{1-\theta}{2\theta}$. Hence
$$
\theta(1+a_0) \ge  \frac{1}{2} + \frac{\theta}{2} > \frac{1}{2}.
$$
Thus \eqref{eq x} has a unique solution $x_0 \in (0,\infty)$.
\qed


\section{Estimate on $\eps^*$}
\label{estimate on eps*}

As mentioned before the largest value of $\eps$ such that \eqref{system00} has a nontrivial solution satisfies
$$
\eps^* \le \eps^0
$$
where
\begin{align*}
\eps^0 = \inf \{ \eps>0 : g_\eps(s) \le 0 \quad \forall s \in [0,1] \} ,
\end{align*}
and
$$
g_\eps(s):=(1-s)f(s)-\eps s .
$$
In the particular case of $f(u) = \chi_{[u>\theta]}$, $0<\theta<1$ we obtain the value $\eps^0 = \frac{1}{\theta}$. 
This estimate can be sharpened.

For this, given $\eps > 0$ let us introduce
$$
G_\eps(b) = \int_0^b g_\eps(s) \,  d s.
$$
For $ 0 < \eps < \eps_0$ define also $a(\eps)$ as the smallest zero of $g_\eps$ in the interval $(0,1)$. We  $b(\eps)$ be the largest $b$ in $[0,1]$ such that $g_\eps$ is positive on $(a(\eps),b(\eps))$.

\begin{prop}
\label{prop eps1}
Assume that $f$ satisfies \eqref{h1} and is continuous. Then \eqref{system00} has no nontrivial solution if $\eps>\eps^1$ where
\begin{align}
\label{def eps1}
\eps^1=\sup \{\eps>0\ \ |\ \ G_\eps(b(\eps))>0\}.
\end{align}
\end{prop}
Before proving the above  proposition, let us recall a result on the scalar problem on a
bounded domain $\Omega \subset \R^N$:
\begin{equation}
\left \{
\begin{aligned}
-\Delta w & = h(w) & \hbox{in } \Omega \\
w & =0 && \hbox{on } \partial \Omega .
\end{aligned}
\right. \label{sca1}
\end{equation}
There  is a vast literature concerning the existence of positive solutions of the above equation.  We will just mention a necessary and sufficient  condition for the existence    of positive solutions.
\begin{thm}
 \label{cn}
(Cl{\'e}ment and  Sweers \cite{CS})
Let $\Omega \subset \R^N$ be a bounded smooth domain and $h \in C^1(\R)$.  Assume there exists $b_1<b_2$ with $b_2>0$ such that $h(b_1)=h(b_2)=0$ and $h>0$ on $(b_1,b_2)$. Then  there exists positive solution $w $ to  \eqref{sca1} with $\max w \in (b_1,b_2)$ if and only if   $\int_{b}^{b_2}h(s)ds> 0$ for all $b \in [0,b_2)$.
\end{thm}

\proof[\bf Proof of Proposition~\ref{prop eps1}]
Fix $ \eps>\eps^1$. From our assumption, we have $G(b(\eps),\eps)< 0$.
Let $\tilde g$ a smooth function such that  $\tilde g_\eps\ge g_\eps$ and $\int_{0}^{\tilde b(\eps)}\tilde g_\eps(s)\,ds<0$.
This is always possible since  $G(b(\eps),\eps)< 0$.
 Observe now that from the necessary and sufficient condition given in theorem \ref{cn}, the following scalar problem
\begin{equation}
  \left \{
   \begin{aligned}
    \Delta w &=-\tilde g_\eps(w) && \hbox{in }  \Omega \\
     w&=0 && \hbox{on } \partial \Omega
   \end{aligned}
  \right.\label{sca2}
\end{equation}
does not have any positive solution for any bounded domain $\Omega$.

Now,  we argue by contradiction, assume there exists a non trivial couple $(u,v)$  solution of the system \eqref{system00}.
An easy computation shows that $u$ is a sub-solution of the problem 
\eqref{sca2} with $\Omega=\R^{n}$.
Take  $\delta >0$ small, and define  $$g_{\eps,\delta}=\tilde g_\eps(s+\delta).$$
So for $\delta$ small enough,(i.e $\delta<\delta_{0}$),we have:
\begin{equation*}
G_{\delta}(\tilde b(\eps)-\delta)=\int_{0}^{\tilde b(\eps)-\delta}g_{\eps,\delta}(s)ds=\int_{\delta}^{\tilde b(\eps)}\tilde g_\eps(s)ds = G(\tilde b(\eps),\eps)-\int_{0}^{\delta}\tilde g_\eps(s)ds\leq 0
\end{equation*}
Therefore $g_{\eps,\delta}$ does not satisfy the necessary existence
condition of Theorem \ref{cn}. Now, take $\Omega=B(0,R)$. Since $u$ goes uniformly to 0 when $|x|$ goes to infinity, we have $\sup_{\partial\Omega} {u}\to 0$ as $R \to +\infty $.
Thus, with $\delta=\sup_{\partial\Omega} {u}<\delta_{0}$, $u_{\delta}=u -\delta$ is a  sub-solution of the  following  scalar problem:
\begin{equation}
  \left \{
   \begin{aligned}
    \Delta w &=-g_{\delta}(w) &&\hbox{in }  B(0,R) \\
    w&=0 && \hbox{on } \partial B(0,R)
   \end{aligned}
  \right.\label{sca3}
\end{equation}
Observe that the constant $\tilde b(\eps)-\delta $ is a super-solution of \eqref{sca3} and
$u_{\delta}<\tilde b(\eps)-\delta$, then we can apply the monotone iterative scheme to
obtain at least one positive solution which  contradicts theorem \ref{cn}.
\qed

\medskip

In the case $f(u) = \chi_{[u>\theta]}$ where $0<\theta<1$ with a similar argument we can show that for $\eps> \frac{(1-\theta)^2}{2}$ the system \eqref{system00}  has no nontrivial solution. Indeed, given any $\sigma>0$ choose a smooth function $\tilde f\ge f$ such that $\tilde f(u) =1$ for $u\ge \theta$ and $\tilde f(u)= 0$ for $u \le \theta-\sigma$. With the same argument as before, \eqref{system00}  has no nontrivial solution for $\eps>\tilde \eps^1$ where $\tilde \eps^1$ is given by \eqref{def eps1} with $f$ replaced by $\tilde f$. A computation then shows that as $\sigma \to 0$, $\tilde \eps^1 \to  \frac{(1-\theta)^2}{2}$.

\bigskip

{\bf Acknowledgment.} J.~Coville warmly thanks Professor Sivashinsky for proposing the formulation of this problem and enlightening discussions.  He is also indebted to Professor Dolbeault for useful conversations on this problem.

J.~Coville was partially supported Tel Aviv University, Universit\'e Paris Dauphine and Max Planck Institute for Mathematics in the Science. 
J.~D\'avila was supported by Fondecyt 1057025, 1090167. Both authors acknowledge the support of Ecos-Conicyt project C05E04.

\end{document}